\documentclass[journal]{IEEEtran}
\usepackage{amssymb,amsmath,amsthm}
\usepackage[usenames, dvipsnames]{xcolor}
\usepackage{tikz}
\usetikzlibrary{arrows,decorations,backgrounds,positioning,fit,petri}
\usepackage{graphicx}
\usepackage{bbm,mathrsfs} 
\usepackage{xspace,braket} 
\usepackage{bm}
\usepackage{cite}
\usepackage{booktabs}
\usepackage[colorlinks, citecolor=blue, linkcolor=black]{hyperref}
\usepackage[capitalize,nameinlink]{cleveref}

\ifCLASSOPTIONcompsoc{}
\usepackage[caption=false,font=normalsize,labelfon
t=sf,textfont=sf, subrefformat=parens, labelformat=parens]{subfig}
\else
\usepackage[caption=false,font=footnotesize,subrefformat=parens,labelformat=parens]{subfig}
\fi

\newtheorem{theorem}{Theorem}
\newtheorem{proposition}{Proposition}
\newtheorem{lemma}{Lemma}
\newtheorem{corollary}{Corollary}

\newtheorem{assumption}{Assumption}

\theoremstyle{remark}
\newtheorem{remark}{Remark}

\theoremstyle{definition}
\newtheorem{definition}{Definition}
\newtheorem{example}{Example}

\newcommand{\VEC}[1]{\bm{#1}}
\newcommand{\mat}[1]{\bm{#1}}
\newcommand{\MC}[1]{\mathcal{#1}}
\newcommand{\MB}[1]{\mathbb{#1}}

\newcommand{\real}{\mathbb{R}}

\newcommand{\tp}{\mathsf{T}}
\newcommand{\norm}[1]{\left\lVert#1\right\rVert}

\DeclareMathOperator*{\minimize}{minimize}
\DeclareMathOperator*{\argmin}{argmin}
\DeclareMathOperator*{\dom}{\mathbf{dom}}
\DeclareMathOperator{\prox}{\mathbf{prox}}
\DeclareMathOperator*{\gra}{\mathbf{gra}\,}
\DeclareMathOperator*{\nul}{\MC{N}}
\DeclareMathOperator*{\ran}{\MC{R}}

\newcommand{\oper}[1]{\ensuremath{\mathsf{#1}}}
\newcommand{\opT}{\oper{T}\,}
\newcommand{\opI}{\oper{Id}\,}
\newcommand{\opR}{\oper{R}\,}
\newcommand{\refl}[1]{\ensuremath{\oper{R}_{\gamma #1}}\,}

\newcommand{\Id}{\mat{I}\xspace}
\newcommand{\vardim}{d}
\newcommand{\ncnt}{l}

\newcommand{\edgeset}{\MC{E}}

\newcommand{\vertexset}{\MC{V}}
\newcommand{\graph}{\mathscr{G}}

\newcommand{\hypergraph}{\mathscr{H}\xspace}

\newcommand{\eigmax}{\lambda_{\textrm{max}}}
\newcommand{\eigmin}{\lambda_{\textrm{min}}}
\newcommand{\eigminnz}{\lambda_{\textrm{min}>0}}

\newcommand{\sigmin}{\sigma_{\mathrm{min}}}

\newcommand{\monotone}{\ensuremath{\alpha}\xspace}
\newcommand{\smooth}{\ensuremath{\beta}\xspace}
\newcommand{\average}{\ensuremath{\theta}\xspace}
\newcommand{\cocoercive}{\ensuremath{\beta}\xspace}

\newcommand{\x}{\VEC{x}\xspace}
\newcommand{\y}{\VEC{y}\xspace}
\newcommand{\z}{\VEC{z}\xspace}
\newcommand{\lag}{\VEC{\nu}}
\newcommand{\s}{\VEC{s}\xspace}

\newcommand{\X}{\VEC{X}\xspace}

\newcommand{\Z}{\VEC{Z}\xspace}

\newcommand{\A}{\mat{A}\xspace}
\newcommand{\B}{\mat{B}\xspace}
\newcommand{\C}{\mat{C}\xspace}
\newcommand{\D}{\mat{D}_v\xspace}
\newcommand{\E}{\mat{D}_e\xspace}

\newcommand{\DCEC}{\D - \C\E^{-1}\C^\tp}
\newcommand{\Lg}{\mat{L}_\graph}

\newcommand{\Q}{\mat{Q}\xspace}

\newcommand{\dscale}{\mat{S}\xspace}

\newcommand{\dr}{\VEC{u}}
\newcommand{\dualform}[1]{\hat{#1}}  

\newcommand{\wtd}[1]{\bar{#1}\xspace}

\newcommand{\wA}{\wtd{\A}}
\newcommand{\wB}{\wtd{\B}}
\newcommand{\wC}{\wtd{\C}}
\newcommand{\wD}{\wtd{\mat{D}}_d}
\newcommand{\wE}{\wtd{\mat{D}}_e}
\newcommand{\wDCEC}{\wD - \wC\wE^{-1}\wC^\tp}

\title{Tight Linear Convergence Rate of ADMM for Decentralized Optimization}
\author{Meng~Ma \quad Bingcong~Li \quad Georgios~B.~Giannakis\\
	Department of ECE and Digital Technology Center, University of Minnesota\\
	Emails: \{maxxx971, lixx5599, georgios\}@umn.edu
	\thanks{This work is supported in part by NSF grants, 10121, 12321, 12345.}
}

\begin{document}

\maketitle

\begin{abstract}
The present paper considers leveraging network topology information to improve the convergence rate of ADMM for decentralized optimization,
where networked nodes work collaboratively to minimize the objective.
Such problems can be solved efficiently using ADMM via decomposing the objective into easier subproblems.
Properly exploiting network topology can significantly improve the algorithm performance.
Hybrid ADMM explores the direction of exploiting node information by taking into account node centrality but fails to utilize edge information.
This paper fills the gap by incorporating both node and edge information and provides a novel convergence rate bound for decentralized ADMM that explicitly depends on network topology.
Such a novel bound is attainable for certain class of problems, thus tight.
The explicit dependence further suggests possible directions to optimal design of edge weights to achieve the best performance.
Numerical experiments show that simple heuristic methods could achieve better performance, and also exhibits robustness to topology changes.
\end{abstract}


\section{Introduction}%
\label{sec:introduction}
Distributed optimization arises frequently in a variety of engineering problems, where $n$ networked nodes work collectively to minimize some objective while maintaining a common decision variable.
The objective is typically an aggregation of nodal local objective that can be private and thus unknown to others.
Collaboration among nodes are allowed when they are connected in the network.
A lot of problems can be cast as decentralized optimization either because required data are dispersed at various locations or the computation is off-loaded to many networked computing machines.
Prominent examples includes  multi-agent coordination, large scale machine learning, distributed tracking and localization, detection and estimation over sensor networks, to name a few~\cite{nedic2009distributed,nedic2010constrained,
	schizas2008consensus,ren2007information,boyd2011distributed,
	giannakis2016decentralized,olfati-saber2007consensus}.

A common strategy to solve distributed optimization problems adopts consensus formulation with a fusion center coordinating all nodes~\cite{sayed2014adaptation}.
This centralized approach is challenged by several issues, such as scalability, load balancing and privacy concerns.
Resources needed for computation and communication at fusion center grow with the network size, which bottlenecks the whole system and makes it vulnerable to single-point failure, hindering its scalability.
Additionally, some private data are not to be shared except for some trusted neighbors.
Therefore, algorithms that only requires \emph{in-network processing}~\cite{rabbat2004distributed} is highly desirable, giving rise to \emph{decentralized optimization}.

Due to its paramount importance, distributed optimization problem has been extensively studied during the past decades.
Many algorithms have been proposed since then.
One particularly popular class is based gradient/subgradient-like methods, including distributed gradient/subgradient descent~\cite{nedic2009distributed}, incremental subgradient~\cite{rabbat2004distributed}, projected subgradient~\cite{sundharram2010distributed,nedic2010constrained}, dual averaging~\cite{duchi2012dual}, gossip~\cite{dimakis2010gossip}.
Typically, a node updates its local decision variable by incorporating information from its neighbors and then taking a gradient/subgradient step based on its local objective.
Depending on the properties of objective function, these methods can be quite slow.
Furthermore, additional assumptions are necessary to ensure convergence, for example, the weight matrix must be doubly stochastic~\cite{nedic2009distributed,sundharram2010distributed,nedic2010constrained}.

In this paper, we consider solving decentralized optimization using the alternating direction method of multipliers (ADMM).
ADMM is particularly appealing for such problems because of its ability to decompose the objective into (easier) subproblems and mild condition for convergence~\cite{bertsekas1989parallel,boyd2011distributed}.
Depending the existence of a central coordinator or not, consensus problems can be solved in either centralized~\cite{boyd2011distributed,bertsekas1989parallel} or decentralized manner~\cite{schizas2008consensus,mateos2009distributed}.
Centralized ADMM suffers from the same issues aforementioned due to the existence of fusion centers, hence the present work focuses on the decentralized case.

In a typical iteration of decentralized ADMM, each node updates its decision variable by aggregating information from its neighbors and minimizing a composite function consisting of local objective and aggregated data.
This process treats each node equally regardless of its topological importance,
which could lead to considerable performance boost.
The hybrid ADMM~\cite{ma2018hybrid} takes advantage of this topology information and designs a local aggregation scheme that has been observed to attain better performance on various networks.
Topological information describes the roles of both nodes and edges, the latter yet to be explored.
This paper proposes the hybrid ADMM over weighted graphs, which not only takes into consideration of relative importance of nodes, but also accounts for edge weights, which, properly utilized, can provide further improvements.

\subsection{Related works}
\textbf{ADMM\@.}
Our work builds on top of the ADMM algorithm.
Since proposed in 1970s~\cite{glowinski1975approximation}, ADMM has been successfully applied to many areas, especially distributed optimization and statistical learning, see~\cite{boyd2011distributed} and references therein.
Decentralized ADMM has widely appreciated merits in applications that requires
\emph{in-network processing}, see~\cite{giannakis2016decentralized} for applications in signal processing.
To take advantage of topology information of nodes, hybrid ADMM~\cite{ma2018hybrid} adopt a local aggregate model where multiple local fusion centers work together each of which coordinates a subset of all nodes.
A closely related work~\cite{iutzeler2016explicit} considers cluster-based ADMM that subsumes decentralized ADMM as special cases, but requires another subprocedure, e.g., gossip algorithm, to compute cluster related quantities.

There are a lot works studying convergence of ADMM since it was proposed~\cite{gabay1976dual}, see also~\cite{boyd2011distributed,bertsekas1989parallel}.
Convergence results for decentralized ADMM did not appear until recently~\cite{hong2017linear,deng2015global,shi2014linear,makhdoumi2017convergence}, due to the difficulty of analyzing the convergence of multi-block ADMM~\cite{lin2015global}.
The relation between convergence rate of ADMM and the topology of underlying graph is explored using lifted Markov chain~\cite{franca2017how}.
Convergence of hybrid ADMM has also been shown~\cite{ma2018hybrid}.

In the setting of distributed averaging, the optimal weights and step size are considered in~\cite{ghadimi2014admm}.
Decentralized ADMM over weighted graphs was also considered~\cite{ling2016weighted}, which offers the flexibility to optimize convergence rate by tuning weights.
The importance of edges was incorporated into hybrid ADMM as well~\cite{ma2018graphaware} and its convergence rate bound was shown to be related to spectral properties of the communication graph.

\textbf{Douglas-Rachford splitting.}
Many convex optimization problems can be analyzed from the perspective of \emph{maximal monotone operators}~\cite{rockafellar1976monotone,bauschke2017convex}.
For example, the problem of minimizing $f$ can be cast as finding the zero of the subdifferntial $ \partial f $, which happens to be a maximal monotone operator provided $f$ is convex, closed and proper.
However, finding zeros of two operators can be difficult, whereas \emph{splitting} methods can be employed to split the problem into easier subproblems~\cite{lions1979splitting,eckstein1989splitting,bauschke2017convex}.
One prominent example is \emph{Douglas-Rachford splitting}~\cite{eckstein1992douglasrachford,douglas1956numerical}.


ADMM was shown to be a special case of a method \emph{Douglas-Rachford splitting}~\cite{gabay1983applications,eckstein1989splitting}.
Therefore, one can prove convergence of ADMM through Douglas-Rachford algorithm.
Recently, a tight linear rate for Douglas-Rachford algorithm has been published~\cite{giselsson2017linear}.
The tight bound for ADMM is deduced from that of Douglas-Rachford algorithm, due to their equivalence.
This connection, however, entails assumptions of the constraints that does not always hold in decentralized case.

\subsection{Contributions}
The present paper propose weighted hybrid ADMM aiming to improve the convergence rate of ADMM for distributed optimization by leveraging topology of the underlying graph.
The proposed WHADMM algorithm not only takes into account the importance of nodes but also considers the significance of edges.
Contributions of the present paper are summarized as following.
\begin{enumerate}
	\item The proposed WHADMM algorithm takes another approach for leveraging network topology via weighted communication graphs, which is different to the HADMM that only considers node importance through local aggregation inside neighborhoods.
	The WHADMM not only extends the HADMM to cope with weighted graphs, but also opens the door of alternative approach to improve convergence speed of ADMM for distributed optimization though tuning edge weights.
	\item Based on the equivalence between ADMM and Douglas-Rachford splitting applied to the dual domain, a novel linear convergence rate bound is provided.
	This rate bound is tight for a certain class of problems, which is the first result in the literature of decentralized optimization using ADMM\@.
	This novel bound explicitly shows the connection between convergence rate and condition of graph Laplacian.
	\item Techniques to find the best weights for certain class of problems are provided. For problems that the technique does not apply, we propose a heuristic method for assigning edge weights.
\end{enumerate}


\subsection{Notation}
We use lower case letters to denote scalar variables and upper case for constants. Vectors and matrices are identified by bold face lower case and bold face upper case letters, respectively.
Let $\A\in\real^{m\times n}$, a matrix with $m$ rows and $n$ columns.
The range of $\A$, denoted by $\MC{R}(\A)$, is the space spanned by the columns of $\A$, i.e., $\ran(\A) = \set{\A\x | \x \in\real^n}$.
Correspondingly, the null space of $\A$, denoted by $\nul(\A)$, is the space of all vectors that yield zero when multiplied by $\A$, i.e., $\nul(\A) = \set{\A\x = \VEC{0}|\x \in \real^n}$.


\section{Problem formulation}%
\label{sec:algorithm}

In this section, we introduce the weighted hybrid ADMM for solving decentralized optimization problems.
We first briefly review hybrid ADMM and how the importance of nodes can be taken into consideration.
Then we proceed to show how to incorporate edge weights.

\subsection{Problem statement}

The distributed optimization problems aim to minimize some \emph{separable} cost function
\[
f(\x) = \sum_{i=1}^n f_i(\x)
\]
by finding a common decision variable $\x\in\real^\vardim$ using $n$ networked computing nodes.
The local cost function $f_i$ typically resides on one computing node with possibly \emph{private} data not to be shared with others~\cite{giannakis2016decentralized} so the overall objective $f$ is not readily available.
Therefore, one fundamental consideration for designing distributed algorithms is to solve this problem efficiently while respecting privacy.

The distributed nature and coupled decision variable makes it difficult for parallel processing.
A common remedy is to decouple the constraint by creating local copies at each node and enforce equality among all nodes.
The resultant problem, which is mathematically equivalent to the original one, can be formulated as
\begin{equation*}\label{pb:originalproblem}
\begin{aligned}
\min_{\{x_i\}} &\quad f(\x) = \sum_{i=1}^n  f_i(\x_i)\\
\text{s.~to} &\quad \x_1 = \x_2 = \cdots = \x_n
\end{aligned}\tag{P1}
\end{equation*}
where $n$ is the total number of nodes and $\x_i\in\real^{\vardim}$ is the local copy of global decision variable $ \x $ at node $i$.
The equality constraints guarantees all copies agree on a single common decision variable, thus obtaining a valid solution of the original problem.
The decoupled problem~\eqref{pb:originalproblem} allows parallel minimization of each local objective while periodically exchanging information among computing nodes to ensure all local copies agree on a common value.

\subsection{Communication graph}
One fundamental constraint on distributed optimization is that
communication between two computing machines may not always be possible, due to limitations such as location, physical connectivity, or privacy concerns.
Such constraints can be conveniently represented by undirected graphs, where each node corresponds to one computing machine and an edge connecting two machines indicates they can exchange information.
An undirected graph is characterized by a tuple $\graph = (\vertexset, \edgeset)$, where $ \vertexset = \{1,2,\ldots,n\} $ is the node set and $ \edgeset \subset \vertexset \times \vertexset $ the edge set.
Edge $(i,j)\in\edgeset$ means node $i$ and $j$ can communicate.
Only undirected graphs without self-loops are considered.

There two types of communication graphs commonly used in distributed optimization: centralized and decentralized.
The centralized communication graphs correspond to the setup where there exists a global coordinator.
All nodes exchange information with the global coordinator only, and no inter-node communication is allowed, which translates to the mathematical formulation $\z = \x_i, \; i = 1, 2,\ldots, n$ upon introducing auxiliary variable $ \z $ for the coordinator~\cite{boyd2011distributed,bertsekas1989parallel}.
Obviously, this requires every node has a communication link to the coordinator that is not always available in practice.
The decentralized communication graphs, however, can be used to model the constraints that only neighborhood communication is allowed.
In this case, no global coordinator exists and nodes can only talk to their immediate neighbors.
Mathematically, this can be captured by $\x_i = \z_{ij}, \: \x_j = \z_{ij}$, where $ \z_{ij}, \z_{ji} $ are auxiliary variables for edge $ (i,j) $~\cite{schizas2008consensus,giannakis2016decentralized}.
Provided that the graph $ \graph $ is connected, we obtain an equivalent formulation of problem~\eqref{pb:formulation}.

\subsection{Hybrid communication graph}

The hybrid communication graph unifies both centralized and decentralized cases, offers insights of distributed optimization problems,
and provides a viable approach to strike a balance between the two aforementioned cases.
The key idea is to introduce \emph{groups}, a generalizations of neighborhood relation, and perform group-based aggregation, instead of aggregating in global or edge-based level.
Nodes residing in the same group can exchange information via group coordinator, thus can be viewed as \emph{group neighbors}.
For example, in centralized communication graph, all nodes belong to the same group and share information through the global coordinator.
While for decentralized case, each edge corresponds to one group consisting of the two connected nodes.
There is no group coordinator in this case since the two nodes can directly talk to each other.

For a given networked system,
numbering groups in arbitrary order and denoting $j$-th group $ \edgeset_j \subset \vertexset $,
a set containing all nodes inside this group.
We say a communication graph is \emph{group connected} if
\begin{equation}
\bigcup_{j=1}^m \edgeset_j = \vertexset
\end{equation}
where $m$ is the total number of groups.
A \emph{path} in a hybrid communication graph is a sequence of nodes, each of which is \emph{group neighbor} of previous one.
Two nodes are said \emph{group-connected} if there is a path from one to the other.

Assigning one auxiliary variable $ \z_j $ to group $j$ and imposing equality constraints between group coordinator and each node, problem~\eqref{pb:originalproblem} is equivalent to
\begin{equation*}\label{pb:formulation}
\begin{aligned}
\minimize_{\{x_i, z_j\}} &\quad  \sum_{i=1}^n f_i(\x_i)\\
\mathrm{subject\: to} &\quad \x_i = \z_j, \; i \in \edgeset_j
\end{aligned}\tag{P2}
\end{equation*}
provided that the graph is \emph{group connected}.

Formulation~\eqref{pb:formulation} is easy to understand, but it also shadows the big picture.
The intrinsic structure of this problem reveals itself once one rearrange the constraints in matrix form.
Let $ \x := [\x_1^\tp, \ldots, \x_n^\tp] \in\real^{n\vardim} $ collect all local copies $ \x_i $ and $ \z := [\z_1^\tp, \ldots, \z_m^\tp] \in \real^{m\vardim} $.
Denote $ \ncnt $ the total number of constraints $\{\x_i = \z_j\}$ for a given problem.
Define $\A\in\real^{\ncnt\times n }$ and $\B\in\real^{\ncnt\times m}$ such that
if the $k$-th constraints is $\x_i = \z_j$, then $A_{ki} = 1, \; B_{kj} = 1$, and $A_{kt} = 0, \; B_{kt} = 0$ for all $ t\neq i, t\neq j$.
With the help of $\A,\B$, constraints in~\eqref{pb:formulation} can be represented compactly by
\begin{equation*}\label{pb:formulation-matrix}
\begin{aligned}
\minimize_{\{x_i, z_j\}} &\quad  \sum_{i=1}^n f_i(\x_i)\\
\mathrm{subject\; to} &\quad
(\A\otimes\Id_\vardim) \x - (\B\otimes\Id_\vardim) \z = \VEC{0}
\end{aligned}\tag{P3}
\end{equation*}
where $\otimes$ denotes Kronecker product.

Two matrices closed related to the hybrid communication graph are \emph{node degree matrix} $\D$, and \emph{group degree matrix} $\E$.
The degree of a node counts the total number of groups it belongs to, and the degree of a group counts the total number of nodes it consists of.
The node degree and group degree matrices are diagonal, i.e., $ \D = \mathrm{diag}(d_1, \ldots, d_n) $, $ \E = \mathrm{diag}(e_1, \ldots, e_m) $, where $d_i$ is the degree of node $i$ and $e_j$ the degree of group $m$.
The structure of $\A$ and $\B$ suggest $\A^\tp\A = \D$ and $\B^\tp\B = \E$~\cite{ma2018hybrid}.

A hybrid communication graph can be completely specified by an \emph{incidence matrix} $\C\in\real^{n\times m}$.
If node $i$ is included in group $j$, then $C_{ij} = 1$; otherwise $C_{ij} = 0$.
The $i$-th row of $\C$ describes all groups node $i$ belongs to, and the $j$-th column of $\C$ describes all nodes group $j$ contains, which immediately yields $d_i = \sum_j C_{ij}$ and $e_j = \sum_i C_{ij}$.
The incidence matrix can be recovered via $\C = \A^\tp\B$~\cite{ma2018hybrid}.


\begin{example}
Consider an distributed optimization problem over a graph shown in Fig.~\ref{fig:example_graph}.
Consider the case $\vardim=1$.
There are 4 nodes in total and 2 groups completely covering all nodes, namely $ \edgeset_1 := \{1, 2\} $ and $ \edgeset_2 := \{2,3,4,5\} $, shown as gray dashed circles.
Apparently, this grouping scheme results in a group-connected graph.
Based on the grouping, the hybrid communication graph can be characterized by linear constraint matrices
\begin{equation*}
\A = 
\begin{bmatrix}
1 & 0 & 0 & 0 & 0\\
0 & 1 & 0 & 0 & 0\\
0 & 1 & 0 & 0 & 0\\
0 & 0 & 1 & 0 & 0\\
0 & 0 & 0 & 1 & 0\\
0 & 0 & 0 & 0 & 1
\end{bmatrix},
\quad
\B = 
\begin{bmatrix}
1 & 0\\
1 & 0\\
0 & 1\\
0 & 1\\
0 & 1\\
0 & 1
\end{bmatrix}.
\end{equation*}
According to the definitions, the node degree matrix, group degree matrix, and incidence matrix are given by
\begin{equation*}
\D =
\begin{bmatrix}
1 & 0 & 0 & 0 & 0\\
0 & 2 & 0 & 0 & 0\\
0 & 0 & 1 & 0 & 0\\
0 & 0 & 0 & 1 & 0\\
0 & 0 & 0 & 0 & 1
\end{bmatrix}
\quad
\E = 
\begin{bmatrix}
2 & 0\\
0 & 4
\end{bmatrix}
\quad
\C=
\begin{bmatrix}
1 & 0\\
1 & 1\\
0 & 1\\
0 & 1\\
0 & 1
\end{bmatrix}.
\end{equation*}
\end{example}

\subsection{Distributed optimization over weighted graphs}



Formulation~\eqref{pb:formulation-matrix} works well for unweighted graphs, thanks to the ability of hybrid communication graph to capitalizing on node importance~\cite{ma2018hybrid}.
But it is incapable to deal with weighted graphs, and more importantly, fails to squeeze more performance out of the underlying network topology.
In the following, we present a formulation that makes it easy to incorporate edge weights.

For simple graphs, an edge is the connection between two nodes, while in hybrid communication graphs, an edge refers to the connection between one node and one group.
Assigning a weight to an edge in hybrid communication graphs boils down to scaling all values transferred via this edge, which translates to applying a scaling matrix to the left of matrices $\A$ and $\B$ in~\eqref{pb:formulation-matrix}, i.e.,
\[
	(\dscale\A\otimes\Id_\vardim) \x - (\dscale{}\B\otimes\Id_\vardim) \z = \VEC{0}
\]
where $\dscale{}\in\real^{\ncnt\times \ncnt}$ is the scaling matrix whose $k$-th diagonal elements $S_{kk}$ is the scale applied to the corresponding edge and $S_{kk}^2 = w_k$ recovers the associated weight of edge $k$.

Let $ \wA = \dscale\A\otimes \Id_\vardim, \wB = \dscale\B\otimes \Id_\vardim $ be the block structured coefficient matrices such that the constraints be written as $\wA\x - \wB\z = \VEC{0}$.
A result similar to~\cite[Lemma 1]{ma2018hybrid} can be shown for weighted graphs.


\begin{proposition}\label{th:matrix-relation}
Let $ \wA $ and $ \wB $ denote the weighted coefficient matrices, then we have
\begin{equation*}
	\wA^\tp \wA = \wD, \qquad \wB^\tp \wB = \wE, \qquad \wA^\tp \wB = \wC
\end{equation*}
where $ \wtd{\C} = \A^\tp\dscale^\tp\dscale\B \otimes \Id_\vardim $, $ \wtd{\D} = \A^\tp\dscale^\tp\dscale\A\otimes \Id_\vardim $, and $ \wtd{\E} = \B^\tp\dscale^\tp\dscale\B\otimes \Id_\vardim$.
\end{proposition}
The proof of \cref{th:matrix-relation} is pretty straightforward.
Simply plugging in the definition of $ \wtd{\A} $ and $ \wtd{\B} $ and using~\cite[Lemma 1]{ma2018hybrid} completes the proof.

Incorporating weight information, the decentralized optimization problem~\eqref{pb:formulation-matrix} can be equivalently formulated as
\begin{equation*}\label{pb:formulation-weighted}
\begin{aligned}
\minimize_{\x_i, \z_j} & \quad \sum_{i=1}^n f_i(\x_i)\\
\mathrm{subject\: to} & \quad \wA\x = \wB\z.
\end{aligned}\tag{P4}
\end{equation*}







\def\len{2cm}
\def\extra{2.5cm}
\begin{figure}[t]
	\centering
	\begin{tikzpicture}[auto,
		vertex/.style={circle, fill=black, inner sep=0mm, 
					   minimum size=.25cm, line width=0pt},
		every edge/.style={semithick, draw}]
	
		\draw[red, thick, rounded corners=.5cm, dashed] 
			(-90:\extra) -- (30:\extra) -- (150:\extra) -- cycle;
		
		\draw[blue, thick, rounded corners, dashed,
			xshift=-\len*1.73*0.75,
			yshift=\len/4,
			rotate=30]
			(-\len/2*1.25, -0.25cm) rectangle (\len/2*1.25, .25cm); 
		
		\node[vertex] (5) at (30:\len) 	[label={[blue]-100:$5$}] 	{};
		\node[vertex] (2) at (150:\len) [label={[blue]-75:$2$}] 	{};
		\node[vertex] (4) at (270:\len) [label={[blue]75:$4$}] 		{};
		\node[vertex] (3) at (0,0)		[label={[blue]above:$3$}] 	{}
			edge (2)
			edge (4)
			edge (5);
		\node[vertex] (1) at (-\len*1.73, 0) [label={[blue]above:$1$}] {}
			edge (2);
	\end{tikzpicture}
	\caption{One possible grouping scheme for a simple graph.
		 Solid black circles represent nodes with their numbering, solid lines indicate connectivity between nodes, and dashed lines identify groups.}%
	\label{fig:example_graph}
\end{figure}
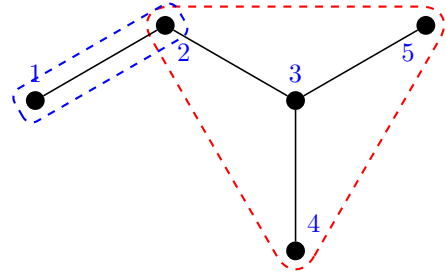

\begin{remark}
Alternatively, one can stack variables $ \{\x_i\} $ and $ \{\z_j\} $ as rows to obtain $ \mat{X} = [\x_1, \x_2, \ldots, \x_n]^\tp \in \real^{n\times \vardim} $ and $ \mat{Z}:=[\z_1, \ldots, \z_m]^\tp \in \real^{m\times \vardim} $.
The constraitns can now be represented by a matrix equality $\A\X - \B\Z = \VEC{0}$.
Notice that we do not need Kronecker product to obtain block structured matrices anymore, effectively reducing the memeory requirements and compuational costs~\cite[Appendix~A]{ma2018hybrid}.
However, using variables as matrices is not that intuitive and further algorithm derivation and convergence proof requires extra efforts to comprehend.
Hence, we keep using variables as vectors for analysis and explanatory purposes.
\end{remark}

\begin{remark}\label{rm:hypergraph}
The hybrid communication graph can also be specified by a hypergraph, which is also a tuple $ \hypergraph = (\edgeset, \vertexset) $, where $ \vertexset:=\{1,2,\ldots,n\} $ is the vertex set and $ \edgeset:=\{\edgeset_i|\edgeset_i \subseteq V\} $ is the hyperedge set.
A hyperedge plays similar roles as a group.
The hypergraph degrades to a simple graph when every hyperedge comprises exactly two nodes, i.e., $|\edgeset_i| = 2$.
\end{remark}

\section{Hybrid ADMM over weighted graphs}
\label{sec:whadmm}

In this section, the weighted hybrid ADMM algorithm for problems over weighted graphs is developed.
In particular, explicit ADMM iterations are given and implementations details featuring in-network processing are also discussed.
Finally, WHADMM is connected to several other algorithms.

\subsection{Weighted hybrid ADMM}
Consider the constrained minimization problem~\eqref{pb:formulation-weighted}.
Attach Lagrange multiplier $ \lag_k \in\real^{\vardim} $ to the $k$-th constraint and collect all multipliers into $ \lag = (\lag_1^\tp, \ldots, \lag_l^\tp)^\tp \in\real^{t\vardim} $.
It turns out that ADMM iterations for solving~\eqref{pb:formulation-weighted} can be carried out using $ \x, \lag $ only with $ \z $ completed eliminated.

\begin{proposition}\label{th:hadmm-iteration}
If $ f $ is convex, $\lag^0 $ satisfies $ \B^\tp \lag^0 = \VEC{0} $, and the Lagrangian of~\eqref{pb:formulation-weighted} admits at least one saddle point, then ADMM iterations for solving~\eqref{pb:formulation-weighted} can be carried out as
\begin{equation}\label{eq:whadmm}
\begin{aligned}
\x^{k+1} &= \argmin_{\x} f(\x) + \frac{\gamma}{2}  \x^\tp\wtd{\D}\x
+ \x^\tp(\y^k - \gamma \wtd{\C}\wtd{\E}^{-1}\wtd{\C}^\tp \x^k)\\
\y^{k+1} &= \y^k + \gamma (\wtd{\D} - \wtd{\C}\wtd{\E}^{-1}\wtd{\C}^\tp)\x^{k+1}
\end{aligned}
\end{equation}
where $ \y^k = \wA^\tp \lag^k $ is a change of variable and $\gamma$ is some algorithmic parameter.
\end{proposition}

\begin{IEEEproof}
Let $f(\x) = \sum_i f_i(\x_i)$.
The augmented Lagrangian of~\eqref{pb:formulation-weighted} is
\begin{equation*}
L(\x,\y,\z) = f(\x) + \lag^\tp(\wA\x - \wB\z) + \frac{\gamma}{2} \norm{\wA\x - \wB\z}_2^2.
\end{equation*}
Minimizing the augmented Lagrangian with respect to each variable, we arrive at
\begin{subequations}
\begin{gather}
-\wA^\tp\lag^k - \gamma(\wA\x^{k+1} - \wB\z^k) \in \partial f(\x^{k+1})\label{eq:min-x}\\
-\wB^\tp\lag^k - \gamma\wB^\tp(\wA\x^{k+1} - \wB\z^{k+1}) = 0\label{eq:min-z}\\
\lag^{k+1} = \lag^k + \gamma(\wA\x^{k+1} - \wB\z^{k+1})\label{eq:min-lag}
\end{gather}
\end{subequations}
where $ \partial f $ is the subdifferential of $f$.

Left multiplying~\eqref{eq:min-lag} by $\wB^\tp$ and adding it to~\eqref{eq:min-z} yields
\begin{equation}\label{eq:Blag}
\wB^\tp \lag^{k+1} = 0, \quad k \geq 0.
\end{equation}
\cref{eq:Blag} holds only for $ \lag^k, \; k \geq 1 $.
With initialization $ \wB^\tp\lag^0 = \VEC{0} $, one can guaranteed $\wB^\tp \lag^k = \VEC{0}$ for all $k\geq 0$.

Using~\eqref{eq:Blag} to eliminate $ \wB^\tp \lag^k $ from~\eqref{eq:min-z},
we obtain closed-form solution of $\z$ update
\begin{equation}\label{eq:update-z}
\z^{k+1} = \wE^{-1} \wC^\tp \x^{k+1}.
\end{equation}
Left multiply~\eqref{eq:min-lag} by $\wA^\tp$ and let $\y=\wA^\tp\lag$, then plug~\eqref{eq:update-z} into~\eqref{eq:min-x} and~\eqref{eq:min-lag}, we arrive at \cref{th:hadmm-iteration}.
\end{IEEEproof}



\begin{remark}
The structure of matrix $\wB$ suggests a nice interpretation of $\wB^\tp\lag^k=\VEC{0}$.
Towards this end, let $\VEC{b} = \B^\tp\lag$.
Then the $j$-th entry $\VEC{b}_j = \sum_i b_{ij} \lag_i = \VEC{0}$; in other words,
Lagrange multipliers associated with group $j$ sum up to zero.
This conclusion holds for all groups and all iterations.
\end{remark}

\begin{remark}\label{rm:y-express-x}
\cref{th:hadmm-iteration} implies $\y^k$ depends on $\{ \x^t \}_{t=1}^k$ and graph topology only, since iteration
\begin{equation}
\y^k = \gamma(\wtd{\D} - \wtd{\C}\wtd{\E}^{-1}\wtd{\C}^\tp) \sum_{t=1}^k \x^t
\end{equation}
relies merely on $\{\x^t\}_{t=1}^k$ and topology related matrices $ \Lg = \wDCEC $, which is a valid Laplacian~\cite{ma2018hybrid}.
In fact, it is exactly the Laplacian of the hypergraph representation of hybrid communication graph (see \cref{rm:hypergraph}), which has been investigated in~\cite{bolla1993spectra}.
\end{remark}

\subsection{Implementation}
The WHADMM algorithm in \cref{th:hadmm-iteration} can be implemented in a decentralized fashion,
thanks to the special structures of $\wtd{\D}, \wtd{\E}$ and $\wtd{\C}$.
To see that, first notice that $\wtd{\D}$ is block diagonal, thus $\wtd{\D}\x$ and $\x^\tp\wtd{\D}\x$ are both fully separable across nodes.
It remains to show $\wtd{\C}\wtd{\E}^{-1}\wtd{\C}^\tp\x$ is also separable.
Towards this end, we first show $\z_j$ can be computed via existing communication graph.
Specifically, for group $j$, $\z_j$ is computed via
\begin{equation}
\z_j = \sum_{i=1}^n \wtd{C}_{ij} \x_i / e_j
= \sum_{i\in\edgeset_j} \wtd{C}_{ij} \x_i / e_j
\end{equation}
which is a weighted sum of all group members.
Let $\VEC{v} = \wtd{\C}\z = \wtd{\C}\wtd{\E}^{-1}\wtd{\C}^\tp\x$.
Then $ \VEC{v}_i $ can be obtained via
\begin{equation}\label{eq:vi}
\VEC{v}_i = \sum_{j=1}^{m} \wtd{C}_{ij} \z_j = 
\sum_{i \in \edgeset_j} \wtd{C}_{ij} \z_j
\end{equation}
which is a weighted sum of $ \z_j $ from all connected group.
Since the computation of $\z_j$ and $\VEC{v}_i$ can be carried out with existing communication graph, we conclude that \cref{th:hadmm-iteration} can be implemented using neighborhood communication only.
In particular, the WHADMM iterations at node $i$ can be accomplished by evaluating
\begin{equation}
\begin{aligned}
\x_i^{k+1} &= \argmin_{\x_i} f_i(\x_i) + \frac{\gamma}{2}d_i \|\x_i\|^2 +  \x_i^\tp(\y_i^k - \gamma\VEC{v}_i^k)\\
\y_i^{k+1} &= \y_i^k + \gamma(d_i \x_i^{k+1} - \VEC{v}_i^{k+1})
\end{aligned}
\end{equation}
where $\VEC{v}_i$ is defined in~\eqref{eq:vi}.

To implement this algorithm in fully decentralized manner, one important questions remains: \emph{who is responsible for computing $ z_j $?}
We have shown that $ z_j $ is maintained per group,
but have not discussed yet which node is actually doing the computation.
It turns out we do not need a dedicate node to fulfill this role, like the global coordinator in the centralized case~\cite{boyd2011distributed}.
Instead by creating groups around some central node, we can create a \emph{virtual fusion center} to perform the tasks.


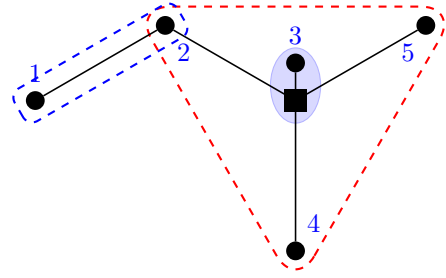
\begin{figure}[t]
	\centering
	\begin{tikzpicture}[auto,
		vertex/.style={circle, fill=black, inner sep=0mm, 
					   minimum size=.25cm, line width=0pt},
		head/.style={rectangle, fill=black, inner sep=0mm, 
					   minimum size=.3cm, line width=0pt},
		every edge/.style={semithick, draw}]
	
		\draw[red, thick, rounded corners=.5cm, dashed]
			(-90:\extra) -- (30:\extra) -- (150:\extra) -- cycle;
		
		\draw[blue, thick, rounded corners, dashed,
			xshift=-\len*1.73*0.75,
			yshift=\len/4,
			rotate=30]
			(-\len/2*1.25, -0.25cm) rectangle (\len/2*1.25, .25cm);

		\filldraw[blue!30, fill opacity=0.5] (90:\len/10) ellipse 
			[x radius=\len/6, y radius=\len/4];
		
		\node[vertex] (5) at (30:\len) 	[label={[blue]-100:$5$}] 	{};
		\node[vertex] (2) at (150:\len) [label={[blue]-75:$2$}] 	{};
		\node[vertex] (4) at (270:\len) [label={[blue]75:$4$}] 		{};
		\node[head]   (0) at (0,0) 		[]							{}
			edge (2)
			edge (4)
			edge (5);
		\node[vertex] (3) at (90:\len/4) [label={[blue]90:$3$}]     {}
		edge (0);
		\node[vertex] (1) at (-\len*1.73, 0) [label={[blue]above:$1$}]   {}
			edge (2);
	\end{tikzpicture}
	\caption{Illustration of \emph{in-network acceleration} on the simple graph from Example~1. The \emph{virtual fusion center} is identified by a black square, while nodes by black circles. The virtual fusion center created inside central node 3 inherits all its connections along with itself.}%
	\label{fig:example_acceleration}
\end{figure}

We term this technique \emph{in-network acceleration} and two ingredients are essential for this recipe to work:
first, groups must be created with a central node connecting all others, limiting the subgraphs induced by groups of diameter two;
and second, the creation of \emph{virtual fusion center} inside the central node which inherits all its connections as well as connects to the central node itself.
The first requirement ensures that the central node is valid group center, so that it can \emph{virtually} act like one.
It should be made clear that the \emph{virtual fusion center} is created conceptually without any physical manifestation, so no dedicated nodes or edges are required.
This technique actually prescribes two roles to the central node, a normal computing node and at the same group center.

For example, in \cref{fig:example_acceleration}, the larger group consists of four nodes and the group center is node 3.
We create a \emph{virtual fusion center}, shown by a black square, and connect it to node 2, 4 and 5 using existing edges, as well as node 3 itself.
The physical node 3 consists of the virtual fusion center and computing node 3.
Thus, WHADMM iterations can be carried out without creating dedicated nodes and edges.




\subsection{Connections and extensions}

WHADMM is closed related to several ADMM-based decentralized optimization algorithms.
When $\dscale = \Id$, it is equivalent to the standard HADMM~\cite{ma2018hybrid}, from which we can further recover several existing algorithms based on the underlying graph topology.
HADMM reduces to centralized ADMM~\cite{boyd2011distributed} when all nodes can be assigned into a single group, i.e., there exists a node that connects all others.
When a group is created for each edge in $ \graph $, HADMM reduces to fully decentralized ADMM when every group corresponds to one simple edge~\cite{giannakis2016decentralized,schizas2008consensus}.
When $ \dscale \neq \Id $, WDADMM can be recovered by creating one group for each simple edge~\cite{ling2016weighted}.
One extension of WHADMM immediately manifests itself through the dependency on edge weights.
For unweighted graphs, WHADMM can apparently benefit from finding the optimal weights for a particular problem.
For weighted graphs, one can simply compensate the value of each weight so that the overall effect is equivalent to the optimal weights.
By incorporating weights into~\eqref{pb:formulation-weighted}, we have transformed weight optimization into the problem of finding optimal preconditioning matrix.
Due to space limit, we leave it for future work.
\section{Convergence analysis}%
\label{sec:convergence_analysis}

In this section, we analyze convergence properties of WHADMM algorithm. 
In particular, we show that WHADMM converges linearly under certain regularity conditions, and provide a new convergence rate bound.
We further demonstrate that this bound is tight for certain class of problems, which is the first results in the literature of decentralized optimization using ADMM\@.


\subsection{Preliminaries}

\textbf{Convex functions.}
The proof relies on some properties of convex functions, which we briefly discuss in the sequel.
A closed, convex and proper function $f$ is $\monotone$-strongly convex with $\monotone >0$ if $f - (\monotone/2)\|\cdot\|^2$ is convex;
$f$ is $\smooth$-smooth with $\smooth>0$ if it is differentiable and $(\smooth/2)\|\cdot\|^2 - f$ is convex.

Conjugate functions play an important role in optimization theories.
The conjugate of a function $ f $ is defined as
\begin{equation}
f^*(\y) = \sup_{\x\in\dom f} \x^\tp \y - f(\x)
\end{equation}
where we use $f^*$ to denote the conjugate of $f$.
If $f$ is strongly convex, then $f^*$ is always differentiable and its gradient is Lipschitz continuous.

\textbf{Monotone operators.}
Let $ \MC{X}, \MC{Y} $ be two nonempty subsets of some Hilbert space and let $ 2^\MC{Y} $ denote the power set of $ \MC{Y} $, the set of all the subsets of $ \MC{Y} $, i.e., $ 2^{\MC{Y}} := \set{\MC{Z}|\MC{Z} \subset \MC{Y}} $.
An operator $ \opT: \MC{X} \mapsto 2^{\MC{Y}} $ is a set-valued relation that maps a point $ \x \in \MC{X} $ to a point $ \MC{Z} \in 2^{\MC{Y}} $.
We use $ \opT \x $ to denote the point $ \x $ is mapped to, and we have $ \opT\x \subset \MC{Y} $.
The graph of an operator $ \opT: \MC{X} \mapsto 2^\MC{Y}$ is defined as
\begin{equation*}
\gra\opT = \Set{(\x, \y) \in \MC{X} \times \MC{Y} | \y \in \opT\x}.
\end{equation*}
We write $ (\x,\y) \in \opT $ in short for $ (\x,\y) \in \gra \opT $ when no ambiguity arises.
If $ \opT \x $ is a singleton, then the mapping $ \opT $ is actually single-valued and can be thought equivalently to a function.
If $ (\x, \y) \in \opT $, then define $ \opT^{-1} $ as the inverse of $ \opT $, which is characterized by its graph $ \gra \opT^{-1} = \set{ (\y, \x) \in 
\MC{Y} \times \MC{X} | \y \in \opT \x } $.

Some special operators that prove useful are briefly discussed here.
We use $ \opI $ to represent the identity operator that maps a point $ \x \in \MC{X} $ to itself, i.e., $ \opI \x = \x $.
Any matrix $ \A \in \MB{R}^{m\times n} $ can be seen as an linear operator that maps a point $ \x \in \real^n $ to $ \A\x \in \real^m $.

The monotone operator theory and convex analysis are closed related~\cite{bauschke2017convex}.
Some properties that bridges these two domains are described in \cref{def:property}.
\begin{definition}\label{def:property}
Let $ \MC{X} $ be a set in some Hilbert space and $ \opT $ an operator $ \opT: \MC{X} \mapsto 2^\MC{X} $.
For any $ (\x, \VEC{u})\in \opT $, $ (\y, \VEC{v}) \in \opT  $, operator $ \opT $ is
\begin{itemize}
	\item \emph{\smooth-Lipschitz} if $ \|\VEC{u} - \VEC{v}\| \leq \smooth \|\x - \y\| $;
	\item \emph{nonexpansive} if $ \opT $ is Lipschitz with $ \smooth \leq 1 $;
	\item \emph{contractive} if $ \opT $ is Lipschitz with $ \smooth < 1 $;
	\item \emph{monotone} if $ \langle \VEC{u} - \VEC{v}, \x - \y \rangle \geq 0 $;
	\item \emph{\monotone-strongly monotone} if $ \langle \VEC{u} - \VEC{v}, \x - \y \rangle \geq \monotone \|\x - \y \|^2 $;
	\item \emph{\cocoercive-cocoercive} if $ \langle \x - \y, \VEC{u} - \VEC{v} \rangle \geq \smooth \| \VEC{u} - \VEC{v} \|^2 $;
	\item \emph{\average-averaged} if $ \opT = (1 - \average) \opI + \average \oper{S} $ for $ \average \in (0, 1) $ and some nonexpansive operator $ \oper{S} $.
\end{itemize}
\end{definition}



One operator that is of paramount importance for analyzing convex functions is the subdifferential operator,
denoted by $ \partial f $, which defines a relation that maps a point $ \x \in \dom f $ to a vector $ \VEC{g} $ satisfying
\[
\gra\partial f = \{ (\x, \VEC{g})| f(\y) \geq f(\x) + 
\VEC{g}^\tp (\y - \x) \}
\]
for all $ \y\in\dom f $.
One can verify that a subdifferntial is always maximal monotone, regardless of the convexity of $f$~\cite{ryu2016primer}.
When $f$ is differentiable, $ \partial f = \nabla f $.

The proximity operator of a function $ f $ is
\begin{equation*}
	\prox_{\gamma f}(\y) := \argmin_{\x\in\dom f} f(\x) + 
	\frac{1}{2} \|\x - \y\|^2
\end{equation*}
and the reflected proximal operator is defined consequently
\begin{equation*}
\oper{R}_{\gamma f}(\y) := 2\prox_{\gamma f}(\y) - \y = (2\prox_{\gamma f} - \opI) \y.
\end{equation*}
Both $\prox_{\gamma f}$ and $\opR_{\gamma f}$ are nonexpansive if $f$ is closed, convex and proper~\cite{bauschke2017convex}.
One benefit of working with reflected proximal is that for a closed, convex and proper function $f$, the fixed points of $\opR_{\gamma f}$ are also the zeros of subdifferential $\partial f$, i.e.,
\begin{equation*}
\x = \opR_{\gamma f} \x \Leftrightarrow \VEC{0} \in \partial f(\x).
\end{equation*}

\textbf{Fixed-point iteration.}
A convex optimization problem often involves minimizing some convex objective function $f$ subject to some constraints.
Finding solutions of such problems is equivalent to searching for a point $\x$ such that $ 0 \in \partial f(\x) $,
which translates the minimization problem to one of finding zero points of subdifferntial operator.
It is well known that if $\opT$ is contractive, then the fixed-point iteration
\begin{equation}\label{eq:fixed-point-iter}
\x^{k+1} = \opT \x^k
\end{equation}
is guaranteed to converge to one of its zero points. 


\subsection{Douglas-Rachford splitting}%
\label{sec:douglas-rachford}

Consider the problem of find zeros of sum of two operators, namely, find a point $\x \in \dom \oper{A}\cup \dom \oper{B}$ such that $0 \in (\oper{A} + \oper{B}) \x$.
Though finding zeros of a single operator $\oper{A}$ or $\oper{B}$ may be straightforward, 
it can be cumbersome to deal with $ \oper{A} + \oper{B} $ as a whole.
Splitting algorithms can decouple the problem into two subproblems, each of which involves only one operator and hence easier to solve.
Among them, Douglas-Rachford splitting (DRS) is one widely used~\cite{eckstein1989splitting,eckstein1992douglasrachford}.

Consider the optimization problem
\begin{equation}\label{eq:min-two-functions}
\minimize_{\x} \quad f(\x) + g(\x)
\end{equation}
where $ f, g $ are both closed, convex and proper.
Let $ \oper{A} = \partial f $ and $ \oper{B} = \partial g $.
Then DRS can be employed to solve this kind of problems.
Let the $ \x^\star $ denote one optimal solution, then it follows that $ \VEC{0} \in  (\partial f + \partial g) (\x^\star) $, which is also one fixed point of $ \opR_{\gamma g}\opR_{\gamma f} $~\cite{eckstein1992douglasrachford}, i.e.,
\begin{equation}\label{eq:fixed-point}
\s^\star = \opR_{\gamma g}\opR_{\gamma f} \s^\star, \quad \x^\star = \prox_{\gamma f} (\s^\star)
\end{equation}
where $ \s \in \dom \oper{R}_{\gamma f} $ is an auxiliary variable to compute $ \x^\star $.

Both $ \opR_{\gamma f} $ and $ \opR_{\gamma g} $ are nonexpansive, so is the composition $ \opR_{\gamma g}\opR_{\gamma f} $.
Solution to~\cref{eq:fixed-point} is computed via fixed-point iteration
\begin{equation}\label{eq:peaceman-rachford}
\s^{k+1} = \opR_{\gamma g}\opR_{\gamma f} \s^k.
\end{equation}
Upon obtaining $\s^\star$, the solution can be retrieved by one proximity step $ \x^\star = \prox_{\gamma f} (\s^\star) $.
This method is known as Peaceman-Rachford splitting (PRS)~\cite{lions1979splitting}.
The problem with PRS is that with convex $ f,\; g $, 
one can only guarantee that $ \opR_{\gamma g}\opR_{\gamma f} $ is nonexpansive, which may not yield a convergent sequence.
PRS is guaranteed to converge when $f$ is strongly convex and $ \opR_{\gamma g}\opR_{\gamma f} $ is contractive.

For a nonexpansive operator, the averaged operator, on the other hand, always generates a convergent sequence. 
If $ \opT $ the averaged PRS operator, i.e., $ \opT = (1 - \average)\opI + \average\opR_{\gamma g}\opR_{\gamma f} $, $\average{}\in (0, 1)$.
then fixed-point iteration of $\opT$ always converges.
The averaged operator $\opT$ is call DRS operator.
DRS reduces to PRS when $\average=1$.
Fixed points of $\opT$ can be found by repeatedly evaluating
\begin{equation}\label{eq:douglas-rachford}
\s^{k+1} = \Big((1 - \average)\opI + \average\opR_{\gamma g}\opR_{\gamma f}\Big) \s^k
\end{equation}
until convergence.
DRS iterations can be carried out equivalently by
\begin{subequations}\label{eq:DR-iteration}
	\begin{align}
	\x^k & = \prox_{\gamma f} \s^k \label{eq:DR-iter-x}\\
	\y^k & = \prox_{\gamma g} (2\x^k - \s^k) \label{eq:DR-iter-y}\\
	\s^{k+1} & = \s^k + 2\average (\y^k - \x^k)
	\end{align}
\end{subequations}
where $ \y, \s $ are intermediate variables facilitating the computation.


If \cref{th:assum-function} holds,
the following results establish tight linear convergence rate of DRS\@.
\begin{assumption}\label{th:assum-function}
	Functions $f$ and $g$ are both closed, convex and proper; moreover, $f$ is $\monotone$-strongly convex and $\smooth$-smooth.
\end{assumption}


\begin{lemma}[{\cite[Theorem~1]{giselsson2017linear}}]\label{th:contraction}
	Suppose \cref{th:assum-function} holds.
	For $\gamma\in(0, \infty)$, $\oper{R}_{\gamma f}$ is $c$-contractive, namely,
	\begin{equation}	
	\| \oper{R}_{\gamma f} \x - \oper{R}_{\gamma f} \y \| \leq c \|\x - \y\|
	\end{equation}
	where
	\begin{equation}\label{eq:contr-factor}
	c = \max\left\{ \frac{\gamma\smooth - 1}{\gamma \smooth  + 1}, \frac{1 - \gamma \monotone}{1 + \gamma \monotone} \right\}
	\end{equation}
	is the contraction factor.
\end{lemma}
For all $\gamma \in (0, +\infty)$, one can verify that $c \in [0, 1)$. 
\Cref{th:contraction} shows that for all positive $\gamma$ and $f$ satisfying \cref{th:contraction}, \refl{f} is contractive.
Together with nonexpansiveness of \refl{g}, the composition $\refl{g}\refl{f}$ is also contractive, ensuring that PRS converges linearly.
\Cref{th:dr-rate} establishes linear convergence rate of DRS and provides explicit rate bound.


\begin{proposition}[{\cite[Proposition 2, 3]{giselsson2017linear}}]\label{th:dr-rate}
	Suppose that \cref{th:assum-function} holds.
	For $\gamma\in(0,\infty)$, $ \average\in(0, 2/(1+c)) $ with $c$ defined in \cref{eq:contr-factor}, the Douglas-Rachford algorithm converges to a fixed point $\bar{\s}\in\mathrm{fix}(R_{\gamma g}R_{\gamma f})$ exponentially at least with rate $\delta$, i.e.,
	\begin{equation}
	\|\s^{k+1} - \bar{\s}\| \leq \delta \|\s^k - \bar{\s}\|
	\end{equation}
	where $ \delta = |1-\average| + \average c $.
	The optimal rate is achieved when setting $ \gamma = 1/\sqrt{\monotone\smooth} $
	\begin{equation}
	\delta^\star = \frac{\sqrt{\smooth/\monotone} + 1 - 2\average}%
	{\sqrt{\smooth/\monotone} + 1}.
	\end{equation}
\end{proposition}



\Cref{th:dr-rate} shows that the rate of convergence of DRS depends on the contraction factor $c$ as well as $\average$.
For a given $f$, $c$ is determined by the strongly convexity and smoothness, see~\eqref{eq:contr-factor}.
But the parameter $\average$ is subject to optimizing.
The closer $\average$ is to one, the smaller $\delta$ is, the faster DRS converges.

The convergence rate bound obtained in \cref{th:dr-rate} is tight, i.e., there exists certain class of problems that this rate can be actually achieved~\cite{giselsson2017linear}.

\subsection{Linear convergence rate of ADMM}%
\label{sec:convergence_rate}

ADMM is known to be equivalent to DRS applied to dual problem~\cite{gabay1983applications}, which offers an alternative approach to study the convergence of ADMM\@.
By exploiting the equivalence, we can extends the tight rate bounds of DRS to obtain a similar result for ADMM\@.
This possibility has been explored in~\cite{giselsson2017linear} and a tight linear convergence rate bound for ADMM is also provided.
However, this result does not hold for all cases since it relies on the condition that $ \A $ has full row rank~\cite[Assumption 2]{giselsson2017linear},
an assumption that is seldom true for the class of decentralized optimization problems~\eqref{pb:formulation-weighted} we are considering in this paper.
In addition, this rate involves only the condition number of $f$ and spectral properties of $\A$, with no explicit connection to the underlying graph.
Therefore, it is difficult to directly transplant the results of DRS to the problem considered here.

To motivate our approach, we first discuss in details why the bound in~\cite{giselsson2017linear} cannot be directly translated to decentralized optimization problem~\eqref{pb:formulation-weighted}.
To simplify the notation, we consider the case $\vardim=1$ in the sequel, which trivially extends to more general cases $\vardim \geq 2$.
With $\vardim = 1$, the constraints of~\eqref{pb:formulation-matrix} reduces to $\A\x - \B\z = \VEC{0}$.

The reason lies in the fact that the strong convexity and smoothness do not carry over to the dual domain when $\A$ does not have full row rank.
In particular, strong convexity of the dual function is shown by bounding $ \lVert \A^\tp\lag \rVert \geq \tau \lVert \lag \rVert $, where $ \tau = \sigmin (\A) $ is the smallest singular,
which only holds when $\sigmin(\A) > 0$.
Without strong convexity of the dual function, they fails to provide any rate bound.
Apart from that, this approach works with matrix $\A$ only, completely overlooking $\B$ the graph topology, hence unable to take full advantage of the problem structures.

We propose a novel convergence rate bound of ADMM that not only holds for surjective $\A$, but also provides an explicit convergence rate when $\A$ is row rank deficient.
When $\A$ has full row rank, this novel rate is equivalent to the existing one~\cite{giselsson2017linear}, which is also tight.
Furthermore, our results forgo the full rank condition and provide an rate bound even when $\A$ does not have full row rank.
This rate has an explicit dependency on graph topology, thus making it possible to design edge weights in order to improve performance.

Our novel bound is achieved by exploiting the structure of problem that has been overlooked by~\cite{giselsson2017linear}.
In particular, we show that it is still possible to bound $ \lVert \A^\tp\lag \rVert $ away from zero even when $\A$ does not have full row rank.
This is achieved by accounting for the fact that dual variables $ \lag $ not only depends on the objective $ f $, but also the graph $ \graph $, see \cref{rm:y-express-x}.
Apart from \cref{th:assum-function}, we also assume the communication graph is connected for otherwise there is no hope to achieve consensus.
\begin{assumption}\label{th:assumption-graph}
	The communication graph is connected, i.e., there is at least one path connecting any two nodes.
\end{assumption}

We first characterize the properties of dual function that is specific to decentralized optimization.
With $\vardim=1$, problem~\eqref{pb:formulation-weighted} can be formulated as
\begin{equation}\label{eq:reformulation}
\begin{aligned}
\max &\quad f(\x) + g(\z) \\
\text{s.~to} &\quad \A\x - \B\z = \VEC{0}
\end{aligned}
\end{equation}
where $ g(\z) = 0 $.
Following the standard procedure, we obtain the dual function
\begin{equation}\label{eq:dual}
\minimize_{\lag} \quad  f^*(-\A^\tp \lag) + g^*(\B^\tp\lag)
\end{equation}
where $f^*$ and $g^*$ are the conjugate functions of $f$ and $g$, respectively.
We choose to write the dual function in minimization form for convenience.

Let $ d(\lag) = d_1(\lag) + d_2(\lag) $, where $d_1(\lag) = f^*(-\A^\tp \lag)$, $d_2(\lag) = g^*(\B^\tp \lag)$.
Since the conjugate of $g$ is $ g^*(\x) = \iota_{\x = \VEC{0}}(\x) $, an indicator function, then $ \oper{R}_{\gamma d_2} = - \opI $ is nonexpansive.
The convergence rate of DRS is determined by the contraction factor of composition $ \opR_{\gamma d_2}\opR_{\gamma d_1} $, which essentially depends solely on $ \opR_{\gamma d_1} $.
Thus, the convergence rate of ADMM can be obtained via characterizing contraction properties of $\opR_{\gamma d_1}$.

\begin{proposition}[{\cite[Proposition~4]{giselsson2017linear}}]\label{th:dual-property-fullrank}
If \cref{th:assum-function} and \cref{th:assumption-graph} hold, and $ \A $ has full row rank,
then $ d_1 $ is $ (\|\A\|^2/\monotone) $-smooth and $ (\tau^2/\smooth) $-strongly convex, where $\tau = \sigmin(\A)$.
\end{proposition}

\cref{th:dual-property-fullrank} claims that smoothness and convexity of dual function rely on spectral properties of $\A$ and objective $f$.
Therefore, for a given function $f$, a well-conditioned $\A$ implies better condition of the dual function.

When $\A$ does not have full row rank, however, using merely nonexpansiveness of $ \opR_{\gamma g^*} $ without accommodating the underlying graph topology leads to no meaningful result, as has been discussed.
This can be circumvented by carefully incorporating problem structure.
To be specific, notice that $ d_2(\lag) = g^*(\B^\tp \lag)$.
Hence, $\dom d = \{\lag | \B^\tp \lag = \VEC{0}\}$, assuming that $ \dom f = \real^\vardim $.
Equivalently, the dual problem can be expressed as
\begin{equation}\label{eq:dual-simple}
\minimize_{\lag\in\dom d} \quad d_1(\lag).
\end{equation}
As discussed in \Cref{sec:whadmm}, the sequences generated by ADMM, which is equivalent to DRS applied to~\eqref{eq:dual}, equip $ \lag $ an expression
\begin{equation}\label{eq:key-relation}
\lag^k = \gamma (\A - \B \E^{-1} \C^\tp) \sum_{t = 1}^k \x^t
\end{equation}
for all $ k \geq 0 $.
The sequence of $ \{\lag^k\} $ in \cref{eq:key-relation} automatically satisfies the constraint $ \B^\tp \lag^k = \VEC{0} $.
With initialization $\lag^0$ such that $ \B^\tp\lag^0 = 0 $, we conclude that
\begin{equation}\label{eq:domain}
\dom d \subset \MC{R}(\A - \B \E^{-1} \C^\tp)
\end{equation}
where $ \MC{R}(\cdot) $ denotes the range of a matrix.

In this specific domain $\dom d$, $\|\A^\tp\lag\|$ can be bounded away from zero, provided that $\lag \neq \VEC{0}$, hence smoothness and strong convexity of $d$ are obtained.
\begin{theorem}\label{th:dual-properties}
If Assumptions~\ref{th:assum-function} and~\ref{th:assumption-graph} hold,
then $ d $ is $\dualform{\smooth}$-smooth and $\dualform{\monotone}$-strongly convex with $\dualform{\smooth} = \|\A\|^2/\monotone$, $ \dualform{\monotone} = \frac{1}{\smooth} \eigminnz^2(\Lg)/\eigmax(\Lg) $, and $ \Lg = \DCEC $ is the Laplacian of the hybrid communication graph.
\end{theorem}

\begin{IEEEproof}
See Appendix~\ref{sec:appendix-dual}.
\end{IEEEproof}
\Cref{th:dual-properties} shows that the dual function $d$ is strongly convex and strongly smooth in its domain, regardless of the row rank of $\A$.
The smoothness depends proportionally to spectral radius of $\A$, which characterizes the maximum node degree of $\graph$, and inverse proportionally to strongly convexity of $f$.
Strong convexity of $d$ is related to smoothness of $f$ and spectral property of graph Laplacian of communication graph.

We have successfully overcome the difficulty brought by singularity of $\A$ and showed strong convexity of the dual even when $\A$ does not have full row rank.
In fact, \cref{th:dual-properties} reduces to \cref{th:dual-property-fullrank} when $\A$ has full row rank.

\begin{corollary}\label{th:dual-property-equivalence}
The strong convexity $\dualform{\monotone}$ defined in \cref{th:dual-properties} is equivalent to that of \cref{th:dual-property-fullrank} when $\A$ has full row rank.
\end{corollary}

\begin{IEEEproof}{}
See Appendix~\ref{sec:proof-equivalence}.
\end{IEEEproof}

With strong convexity and smoothness of the dual function, we are ready to show the contraction factor of $\opR_{\gamma d_1}$ and linear convergence rate of ADMM follows.

\begin{lemma}\label{th:dual-contraction}
If Assumptions~\eqref{th:assum-function} and~\eqref{th:assumption-graph} hold, and $ \gamma \in (0, +\infty) $,
then $\opR_{\gamma d_1}$ is $ \dualform{c} $-contractive, i.e.,
\begin{equation*}	
	\| \oper{R}_{\gamma f} \x - \oper{R}_{\gamma f} \y \| \leq c \|\x - \y\|
\end{equation*}
where the contraction factor
\begin{equation}
\dualform{c} = \max \left\{ \frac{\gamma \dualform{\smooth} - 1}{\gamma \dualform{\smooth} + 1}, \frac{1 - \gamma\dualform{\monotone}}{1 + \gamma\dualform{\monotone}} \right\}.
\end{equation}
\end{lemma}
\begin{IEEEproof}
Applying \cref{th:contraction} to the dual function and using \cref{th:dual-properties} yield the result.
\end{IEEEproof}


\begin{theorem}\label{th:main-result}
If Assumptions~\eqref{th:assum-function} and~\eqref{th:assumption-graph} hold,
and $\gamma\in(0, +\infty)$, $\average \in (0, 2/(1 + \dualform{c})$,
then DRS for solving problem~\eqref{pb:formulation-weighted} converges linearly towards some fixed point at least with rate $ \hat\delta = |1-\average| + \average \dualform c $, i.e.,
\begin{equation}
\|\dr^{k+1} - \bar{\dr}\| \leq \hat\delta \|\dr^k - \bar{\dr}\|
\end{equation}
where $\dr$ is some auxiliary variable and
$ \bar \dr \in \mathrm{fix}(\oper{R}_{\gamma d_2} \oper{R}_{\gamma d_1}) $ is some fixed point.
The dual variable can be recovered via $ \lag^k = \prox_{\gamma d_1} \dr^k $.
\end{theorem}

\begin{IEEEproof}
See Appendix~\ref{sec:proof-main-result}.
\end{IEEEproof}

\Cref{th:main-result} establishes the linear convergence rate of DRS applied to the dual problem and provides explicit form of this rate.
One pitfall is that the convergence rate is derived for $\VEC{u}$ instead of Lagrange multipliers $\lag$.
To establish the convergence rate of ADMM iterations,
denote $ \lag^\star $ one optimal solution, then sequence $ \{\lag^k\} $ satisfies
\begin{equation}
\|\lag^k - \lag^\star \| \leq (\lvert 1 - \average\rvert + \dualform{c} \average)^k
\frac{1}{1 + \gamma\sigma} \|\dr^0 - \bar{\dr}\|
\end{equation}
which proves linear convergence of ADMM\@.

Similar to \cref{th:dr-rate}, there exists some optimal $\gamma$ that leads to optimal theoretical rate.
\begin{corollary}
The parameter $\gamma$ that maximizes the rate bound is $ \gamma = 1/\sqrt{\dualform{\smooth} \dualform{\monotone}} $, with the optimal rate bound given by
\begin{equation}
\dualform\delta^\star = \frac{\sqrt{\dualform\smooth/\dualform\monotone} + 1 - 2\theta}{\sqrt{\dualform\smooth/\dualform\monotone} + 1}.
\end{equation}
\end{corollary}



\begin{remark}
Strictly speaking, ADMM corresponds to DRS with $ \average = 1/2 $ applied to the dual problem. 
When $ \average = 1 $, DRS reduces the PRS algorithm, which enjoys fastest convergence rate when it converges.
When $ \average > 1 $, it is equivalent over relaxed ADMM~\cite{eckstein1992douglasrachford}.
\end{remark}

\begin{remark}
It's not surprising $ \delta^\star $ depends on the cost function and graph topology.
This observation suggests it is possible to improve WHADMM by optimizing graph condition number via edge weights tuning.
\end{remark}

\begin{remark}
	The proof of \cref{th:dual-property-equivalence} shows that when the whole network belongs to one single group, the Laplacian is equivalent to a network of the same nodes with each pair of nodes connected.
	This observation provides insights into this algorithm and sheds light on the performance improvement brought by groups.
	The Laplacian of a single-group graph is the same as a complete graph with the same nodes, suggesting that by creating a group one essentially transform this network into a complete graph.
	If we measure the connectivity of a graph by the smallest nonzero eigenvalue of its Laplacian, then the well known results from spectral graph theory indicate that complete graphs are better, or at least equally, connected than any other graph.
\end{remark}

\subsection{Tightness of the bound}
It has been shown that the linear convergence rate in \cref{th:main-result} is tight when $\A$ has full row rank, i.e., there exist a certain class of problems satisfying the assumption can actually achieve this rate of convergence, see~\cite{giselsson2017linear} for proof and examples.
\Cref{th:dual-property-equivalence} shows the convergence rate bound in~\cref{th:main-result} is equivalent to that of~\cite{giselsson2017linear} when all the nodes can be assigned to one group, which immediately implies our bound is also tight.
For general graphs that do not fit into one group, however, this bound typically cannot be achieved.
\section{Numerical experiments}%
\label{sec:numerical_experiments}

\begin{figure*}[t]
\centering
\subfloat[]{\includegraphics[width=0.3\textwidth]{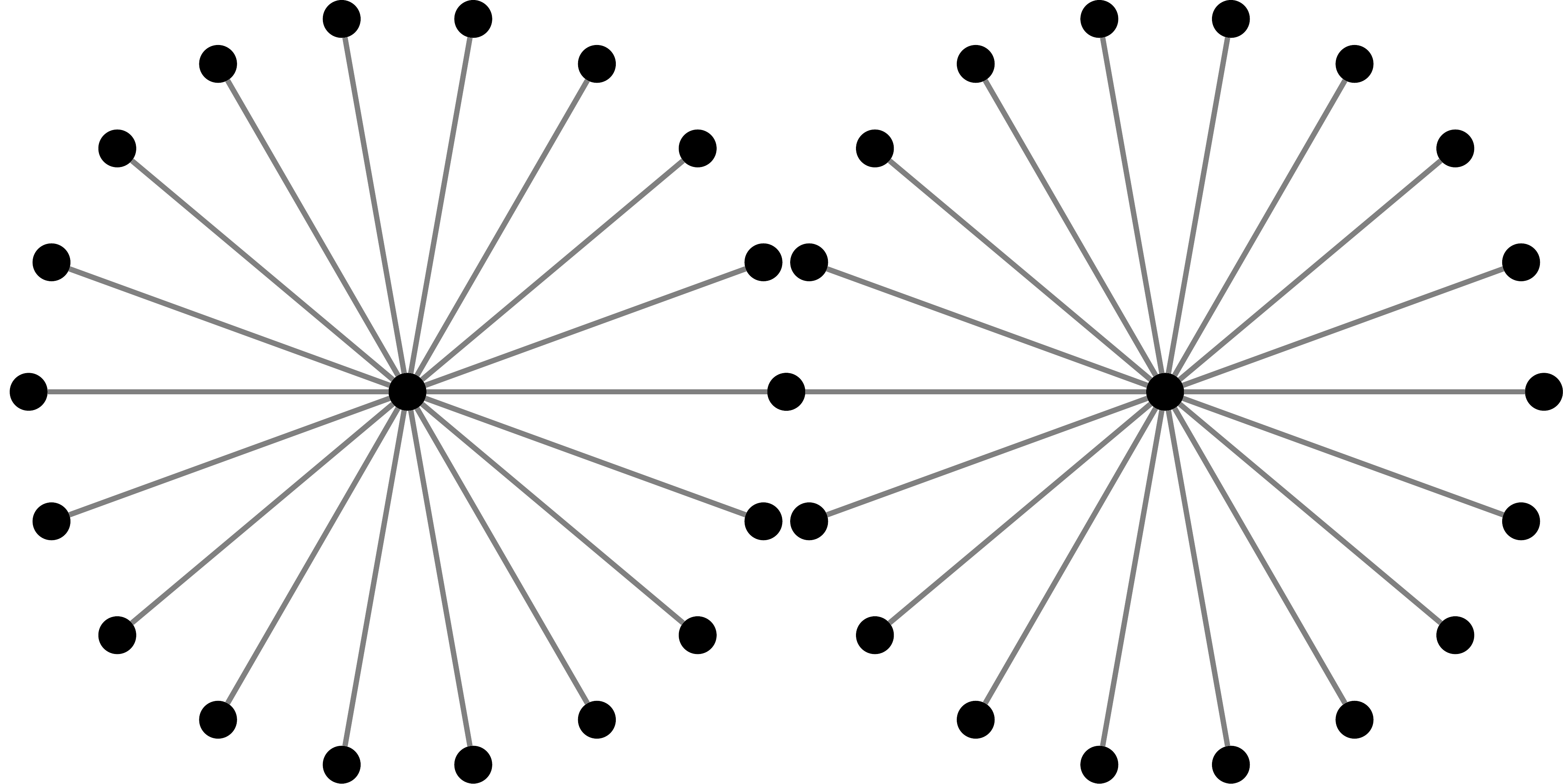}\label{fig:twocluster1}}
\hfil
\subfloat[]{\includegraphics[width=0.37\textwidth]{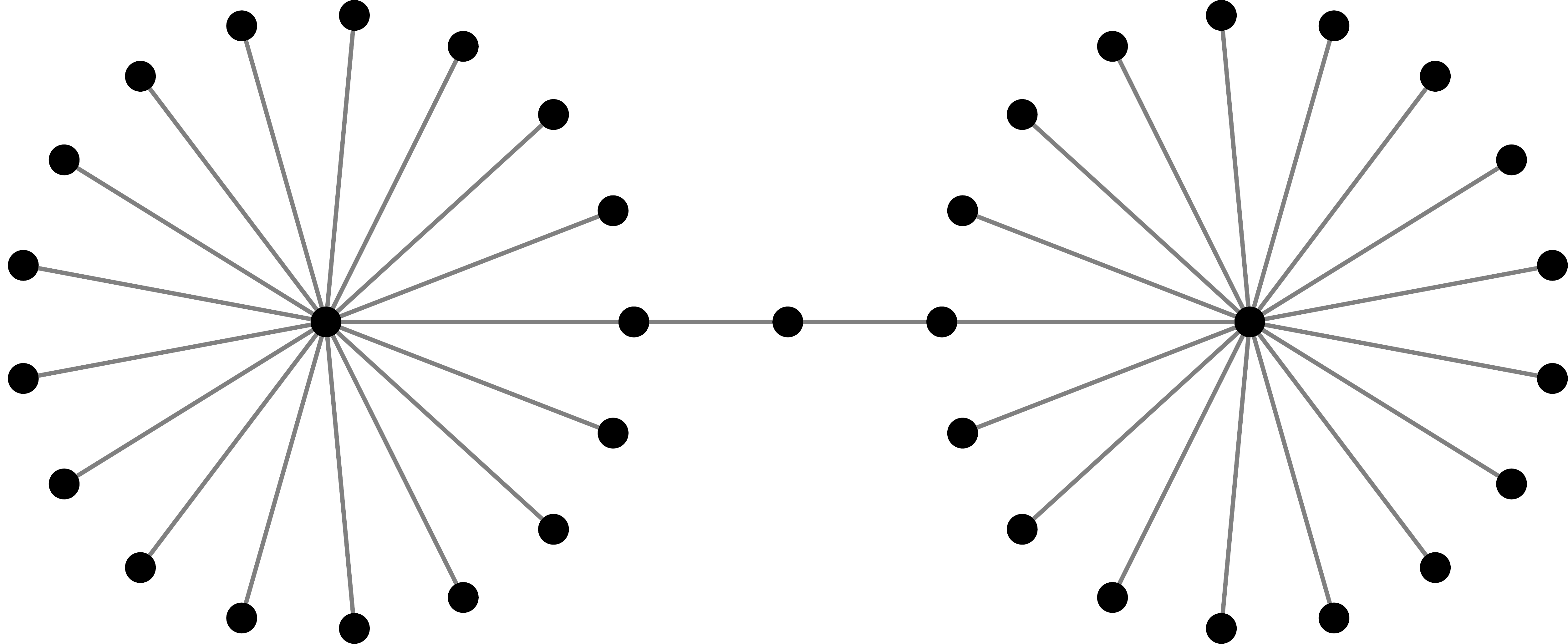}\label{fig:two_cluster2}}
\caption{Two types of graphs considered in the numerical experiments. Both graphs consist of the same number of nodes arranged in two clusters. In the left, two centers of clusters are connected by a path of length 2, while in the right graph 4.}
\label{fig:twocluster}
\end{figure*}

In this section, we perform numerical experiments to demonstrate the linear convergence rate of weighted hybrid consensus ADMM for decentralized optimization, 
and compare with hybrid ADMM to show the possible performance improvement brought by carefully designed weights.
We also compare with ADMM without grouping and local fusion centers and its weighted version, to showcase the benefits of leveraging the importance of nodes.

\textbf{Problem setup.}
The problem tested is weighted network average.
For a graph of $n$ nodes, suppose each node can make a observation.
The weighted network average problem aims to find the weighted average of all nodal observations.
Denote $\VEC{b}_i\in\real^\vardim $ the nodal observation, $\vardim = 3$.
Elements of $ \VEC{b}_i $ are randomly drawn from the interval $ [0, n] $ according to uniform distribution.
Let $q_i$ be the weight of node $i$, which is drawn independently from the interval $[1, \kappa]$, where $\kappa$ controls the condition number.
In the following tests, we let $\kappa=5$.
Theoretically, the weighted average is given by $\sum_{i=1}^n q_i \VEC{b}_i/\sum_{i=1}^n q_i$.

The solution can be found by solving an optimization problem.
Define the local objective function at node $i$ as $ f_i(\x) = q_i \|\x - \VEC{b}_i\|^2 $ and $ x \in \real^\vardim $ is our estimate of the weighted average.
The overall goal is to find the weighted average of all observations $\VEC{b}_i$ by minimizing
$\sum_{i=1}^{n} f_i(\x) = \sum_{i=1}^n q_i\|\x - \VEC{b}_i\|^2 $.
Collecting $\{\VEC{b}_i\}_{i=1}^n$ into rows of $\mat{B} \in \real^{n\times \vardim} $ and $\{q_i\}_{i=1}^n$ as diagonal elements of $\Q\in\real^{n\times n}$, we can write the overall objective in a compact form
\begin{equation}\label{eq:wna-orignial}
\minimize \; \mathrm{Tr}
\left((\VEC{1}\x^\tp - \mat{B})^\tp \Q (\VEC{1}\x^\tp - \mat{B})\right)
\end{equation}
where $\mathrm{Tr}(\cdot)$ computes the trace of a matrix.
By setting the gradient to zero,  
one can easily obtain that $ \x^\star = \mat{B}^\tp \Q \VEC{1}/\mathrm{Tr}(\Q) $ minimizes problem~\eqref{eq:wna-orignial}, which is exactly the weighted average of all $\VEC{b}_i$.

To solve this problem in decentralized manner, create a local copy $ \x_i $ for each node to decouple the dependence on global decision variable $\x$, and
problem~\eqref{eq:wna-orignial} can be reformulated as
\begin{equation}
\begin{aligned}
\minimize & \quad \mathrm{Tr}
\left((\mat{X} - \mat{B})^\tp \Q (\mat{X} - \mat{B})\right)\\
\text{subject to} & \quad \A\X - \B \Z = \VEC{0}.
\end{aligned}
\end{equation}
where $ \X $ collects all local copies $ \x_i $ as rows and $\A, \B$ are created according to the underlying communication graph.

\textbf{Algorithm setup.}
We solve this problem using decentralized ADMM\@.
In particular, we compare the performance of four algorithms, namely, HADMM, WHADMM, DADMM and WDADMM\@.
Both HADMM and WHADMM use \emph{degree centrality} as the centrality measurement of node importance, based on which groups are created.
Other centrality measurements are also possible but not tested here since our focus is to demonstrate the effects of taking into account edge importance.

Due to space limit, we do not explore systematic ways to characterize weight importance, thus we propose an heuristic method that employs edge betweenness centrality~\cite{brandes2001faster} as edge weights.
Intuitively, edge betweenness centrality measures the relative importance of an edge by counting the number of shortest paths between any two nodes that passes through the very edge.
Critical edges, for example the path connecting two components, have larger betweenness centrality, without which the graph becomes disconnected.
Such edges happen to be critical in decentralized optimization as well.
Therefore, we propose to use \emph{normalized edge betweenness centrality} as a heuristic approach for designing edge weights.
Using the \emph{in-network acceleration} technique, each node takes advantage of the existing connection to group center and reuses it as the edge to \emph{virtual fusion center}.
For example, in Fig.~\ref{fig:example_acceleration}, node 2 has two connections, one to node 1 and one to the virtual fusion center (square).
As a result, node 2 can reuse use the betweenness centrality of edge $(2, 3)$ as the weight of its connection to virtual fusion center. 
The group center, node 3, is also connected the virtual fusion center but has no physical correspondence.
We use unit weight for such cases.

To better illustrate the effect of edge weights on algorithm performance, we tests all four algorithms on graphs with two clusters.
Specially, we consider two types of graphs, shown in Fig.~\ref{fig:twocluster}.
Two types of graphs share the same number of nodes arranged in two clusters.
The centers of two clusters are connected by a path.
The only difference between two types is the length of the path, which leads to different connectivity.

To measure the progress while solving this problem, we use the \emph{relative error} $ \|\x - \x^\star\|/\|\x^\star\| $ to indicate how close our solution is to the optima at certain iteration.
We record the \emph{relative error} and corresponding iteration number for all algorithms.
Though~\cref{th:dr-rate} provides the ``best'' $\gamma$ in theory, we observe that it is not always the best in practice, especially when $\A$ is row rank deficient.
So instead, we choose the best $\gamma$ using grid search for each algorithm and graph individually.

\begin{figure*}[t]
\centering
\subfloat[]{\includegraphics[width=0.3\textwidth]{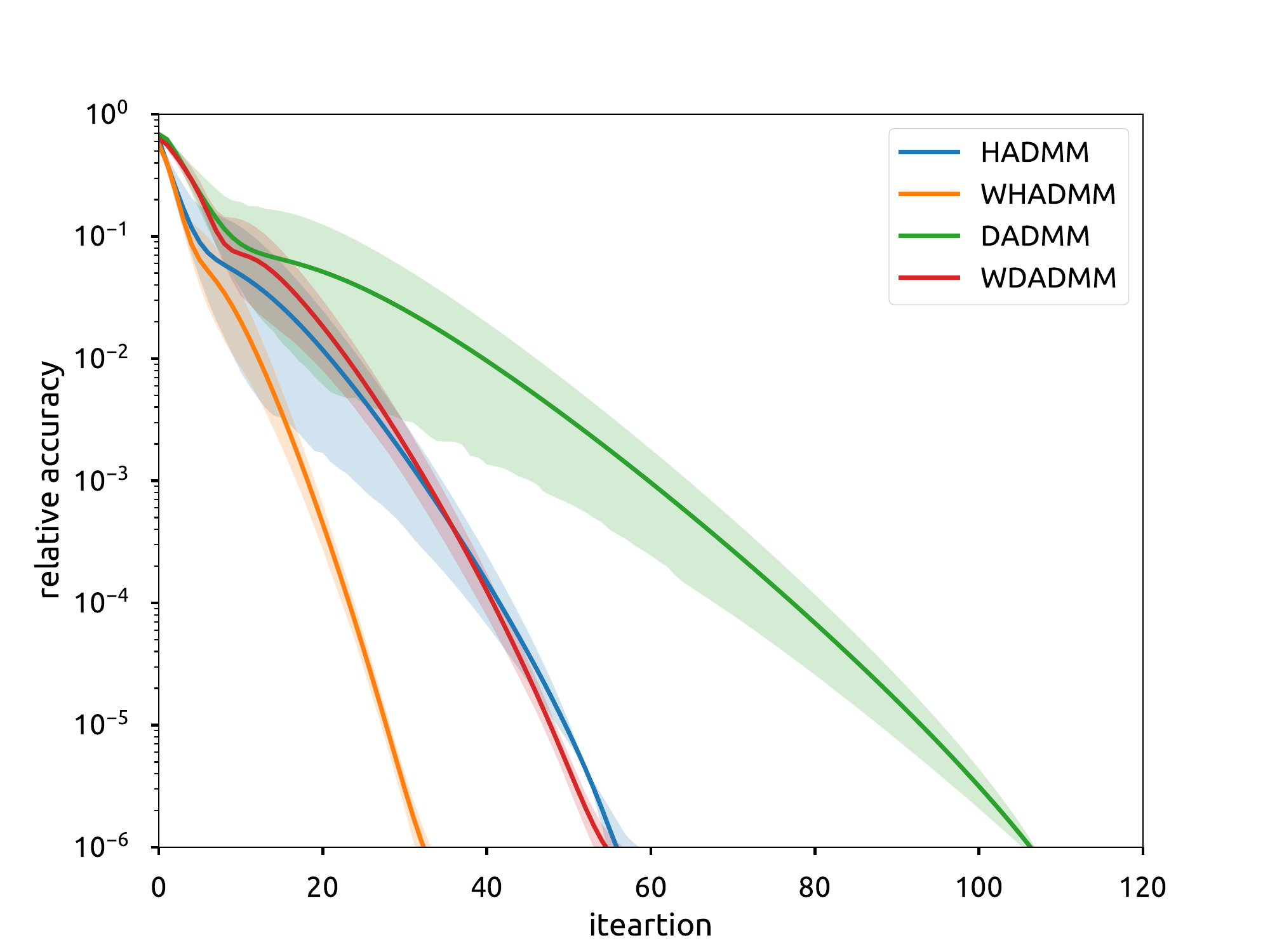}}
\hfil
\subfloat[]{\includegraphics[width=0.3\textwidth]{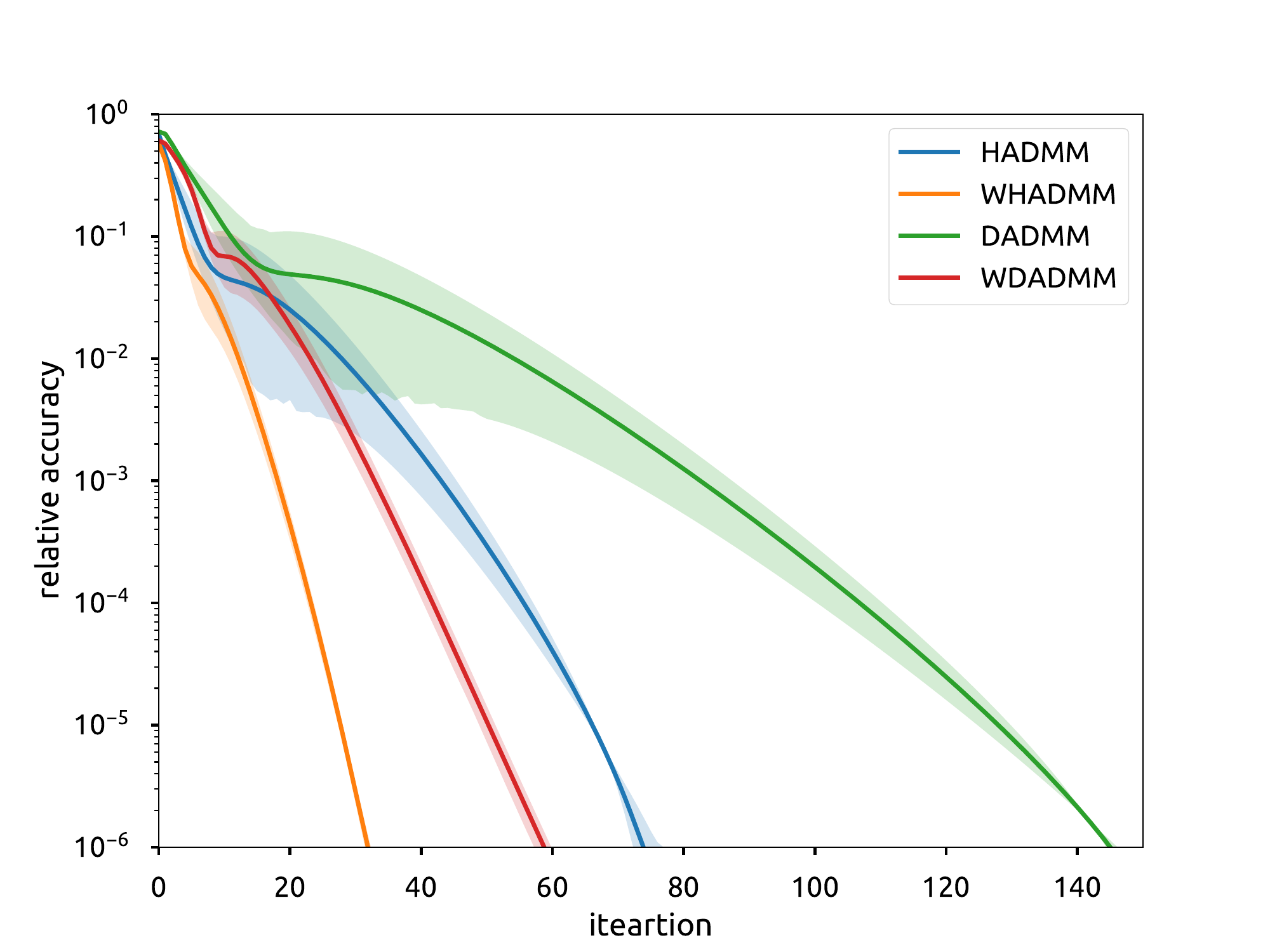}}
\hfil
\subfloat[]{\includegraphics[width=0.3\textwidth]{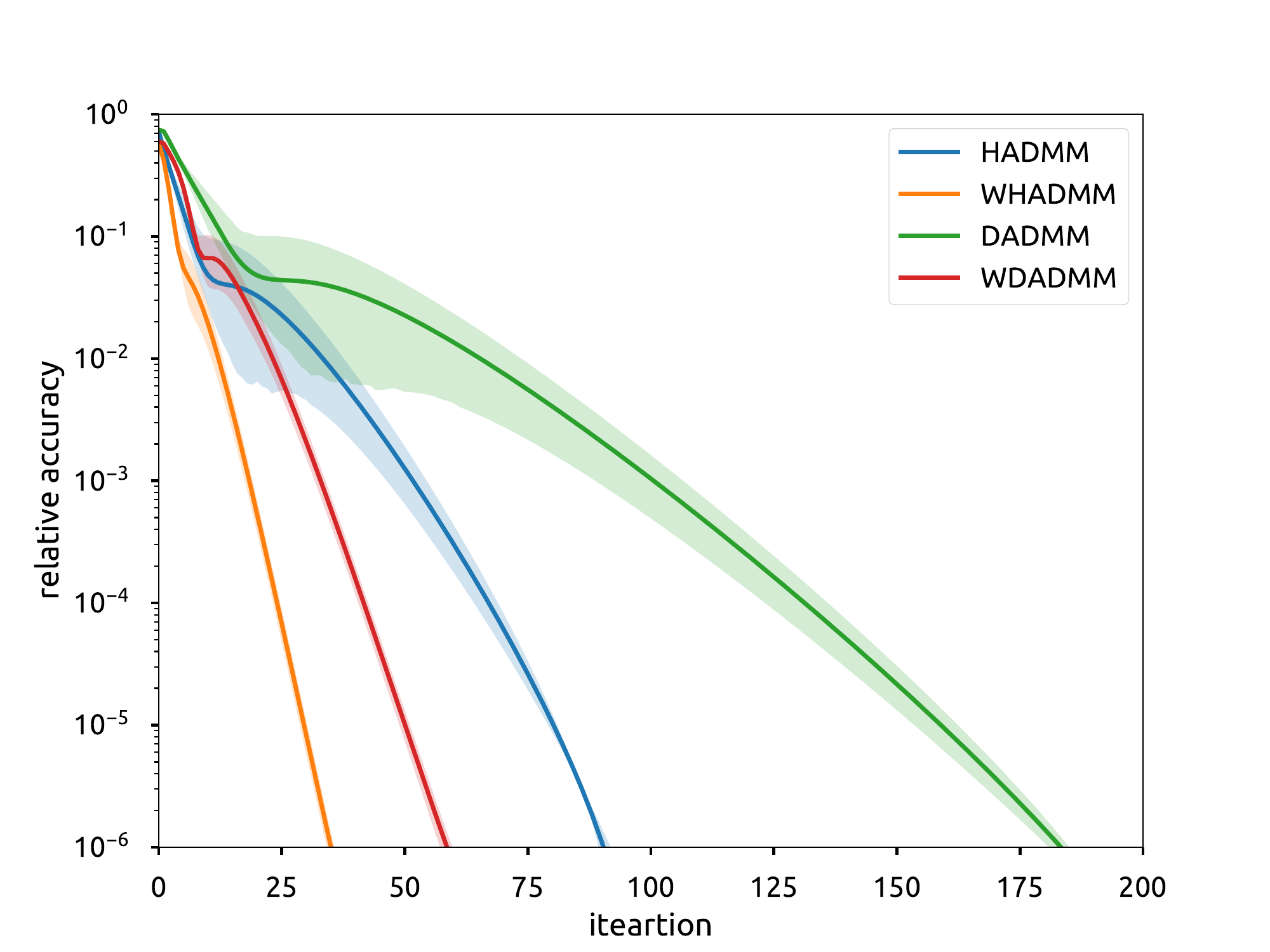}}\\
\subfloat[]{\includegraphics[width=0.3\textwidth]{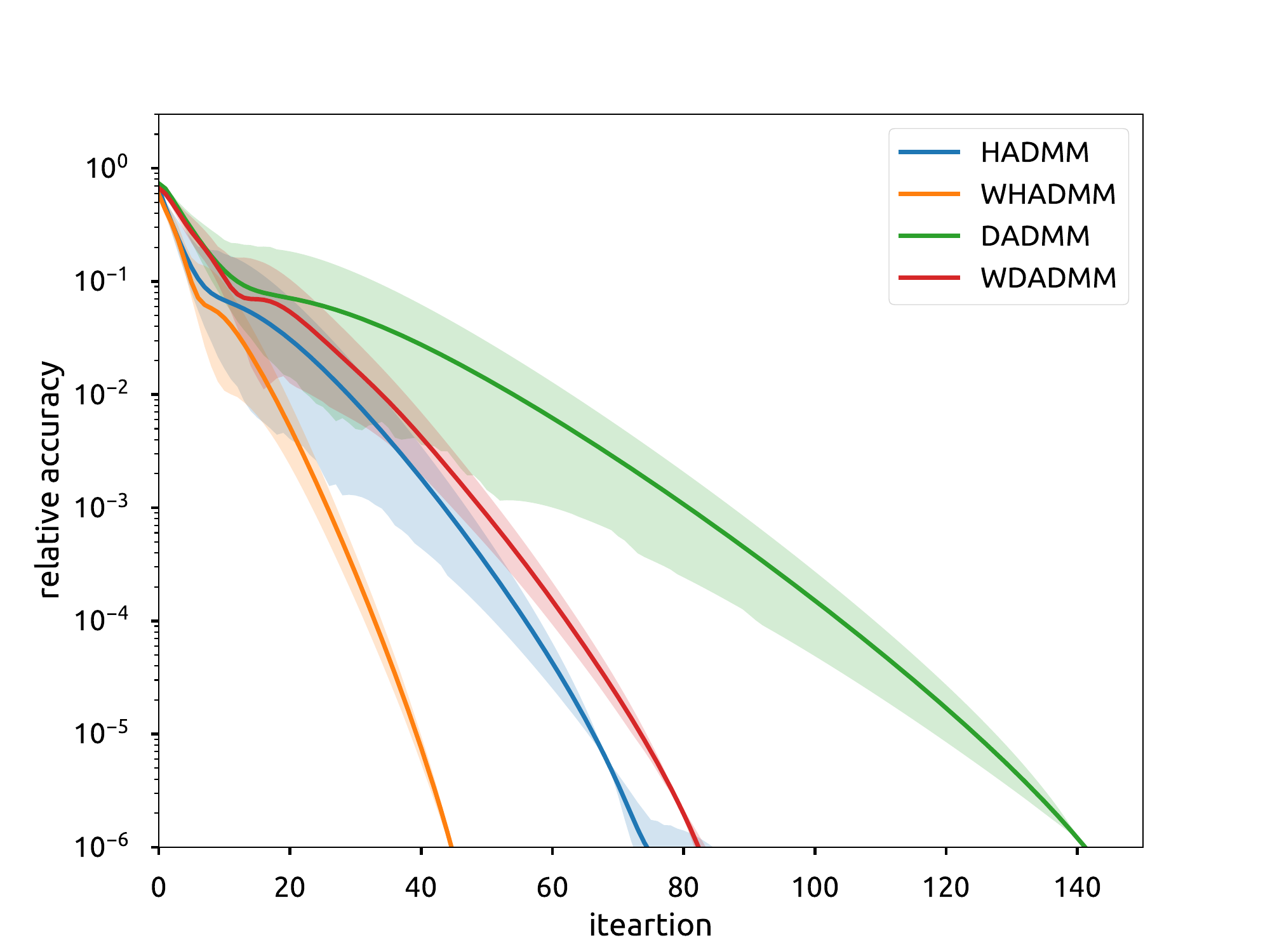}}
\hfil
\subfloat[]{\includegraphics[width=0.3\textwidth]{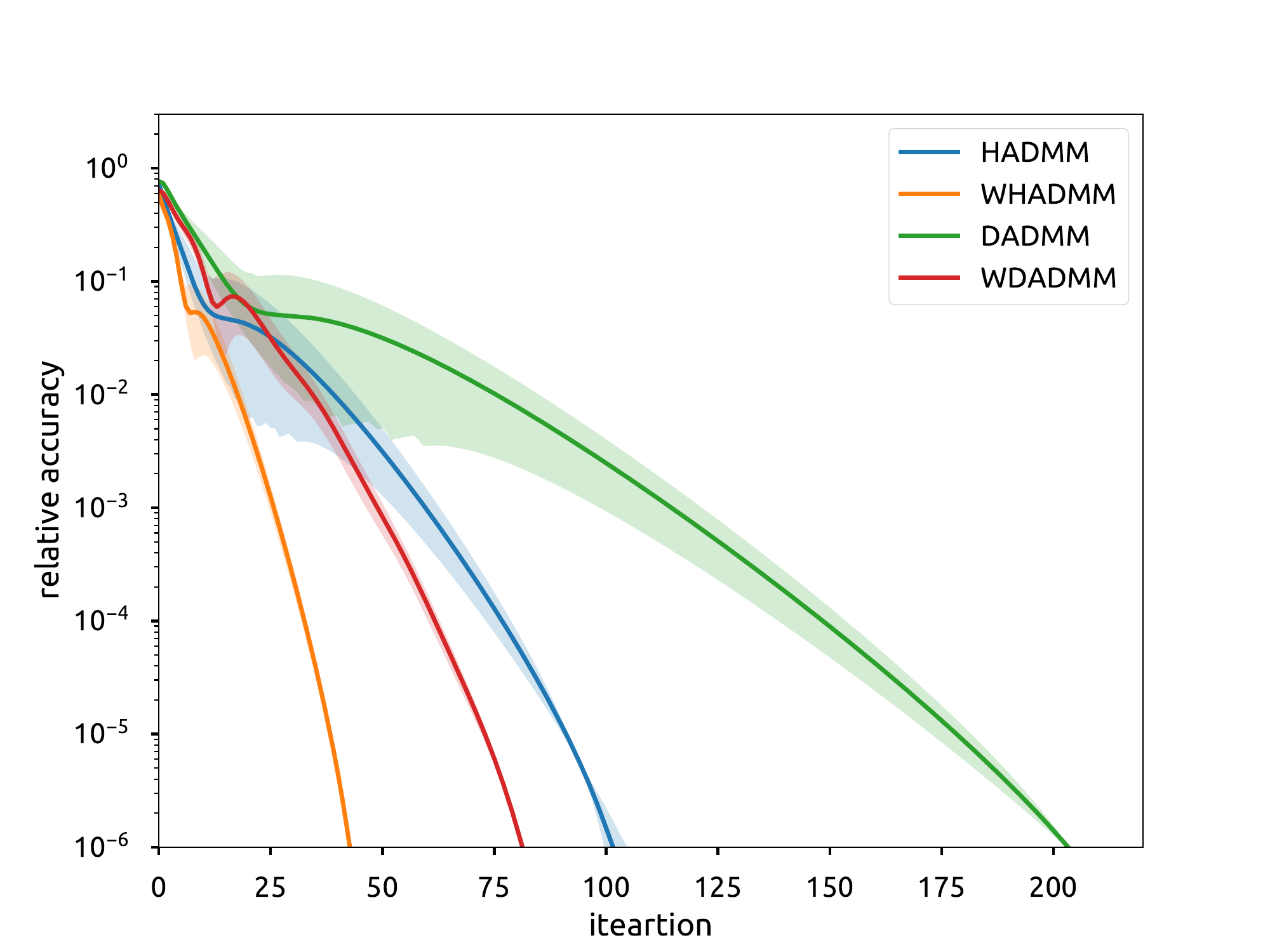}}
\hfil
\subfloat[]{\includegraphics[width=0.3\textwidth]{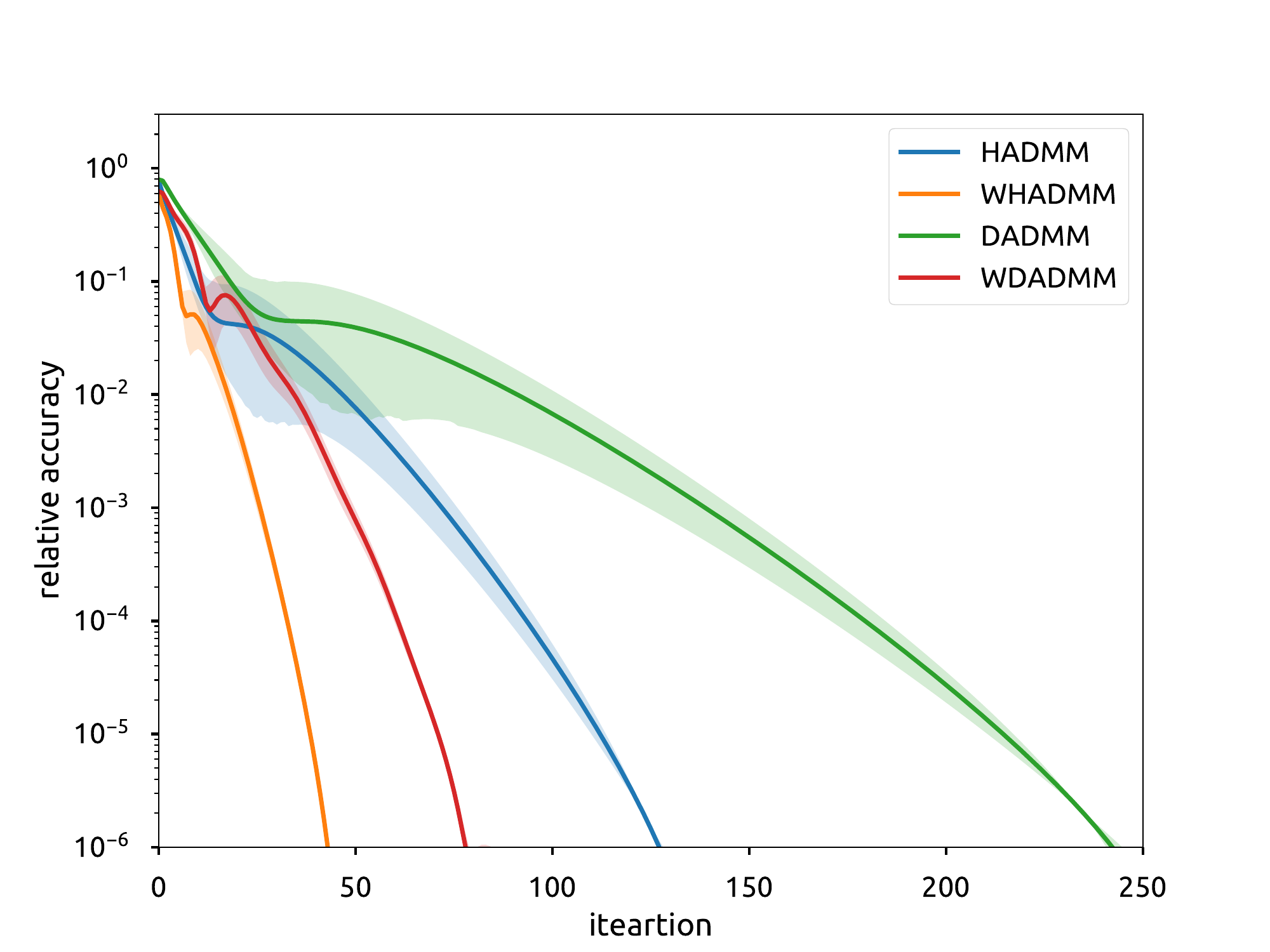}}
\caption{Convergence of hybrid ADMM over graphs with two clusters.
The two graphs contain the same number of nodes.
The only difference is two cluster centers are joined by a longer path in the second graph.
For comparison, performance of DADMM and WDADMM are also plotted.}%
\label{fig:performance-twocluster}
\end{figure*}


\textbf{Results.}
We plot the relative error vs\@. iteration number for all algorithms and all graphs, and the results can be found in Fig.~\ref{fig:performance-twocluster}.
The top row are results on the first type of graphs and the bottom row the second type of graphs.
The first, second and third columns corresponds to the performance on graphs of 21, 41 and 61 nodes, respectively.
Due the randomness, the performance tends to fluctuate across different runs.
To diminish the effect of randomness, we employ Monte Carlo method and independently run each experiment 500 times.
Based on the results, we plot the average in solid line as a representation of this algorithm's performance, and shade the area defined by 5-th and 95-th percentile of all Monte Carlo runs in order to compare their robustness to random observations.

With a simple glance of all figures, we can quickly make two observations.
First, hybrid methods, including HADMM and WHADMM, outperform their conventional counterparts, including DADMM and WDADMM, respectively.
Secondly, weighted methods, including WHADMM and WDADMM, outperform their unweighted counterparts, including HADMM and DADMM\@.
These two observations corroborate the two aspects of levering network topology in decentralized optimization, namely, importance of nodes, which is taken care of hybrid updates, and importance of edges, which is taken care of by weighted updates.
Comparing HADMM and WDADMM with DADMM, we conclude that each of the two aspects brings some performance improvement, and comparing WHADMM with the others we see that combining both offers the most enhancement.

Comparing the columns of figures, we observe a trend on both types of graphs that unweighted algorithms suffer more as the graph size increases, while weighted ones remains stable.
The number of iterations needed to reach a $10^{-6}$ solution for DADMM and HADMM almost double when graph size increased from 21 to 61.
But for WHADMM and WDADMM the number is almost the same.
The same pattern can be observed in both types of graphs.
This clearly demonstrates the value of edge weights in decentralized optimization.

Compare the rows of figures, we see that all algorithms suffer from the worse connectivity of second type graphs.
But still, weighted algorithms are much more stable and suffer far less than unweighted ones.

One may want to compare HADMM with WDADMM to see how the two aspects of leveraging network topology compare.
Unfortunately, there is no definite answer.
For smaller graphs, HADMM may converges faster than WDADMM, but for larger graphs the opposite is true.
So we cannot really say anything regarding this, but at least both are important.

Finally, we check the robustness to random observations.
The same pattern shows up again, namely, unweighted algorithms are much more fragile than weighted ones, and conventional algorithms are more vulnerable then hybrid ones.

\section{Conclusions and discussion}

This paper introduced WHADMM that is able to deal with weighted graphs and has the advantages of exploiting topology information to achieve better performance.
By maximizing the convergence rate, the optimal weights obtained carry a meaningful interpretation, in the sense that critical edges are regarded as more important, aligned with intuition.
Future works include exploring dropping edges by assigning zero weights in order to reduce communication overhead, and extending beyond ADMM to other optimization methods.

\appendices%

\section{Proof of Theorem~\ref{th:dual-properties}}%
\label{sec:appendix-dual}
If $f$ is $\monotone$-strongly convex and $\smooth$-smooth, then the conjugate $f^*$ is $1/\smooth$-strongly convex and $1/\monotone$-smooth, which implies $\nabla f^*$ is $1/\monotone$-Lipschitz continuous
The Lipschitz continuity of $ \nabla f^* $ leads to
\begin{equation}
\begin{aligned}
&\; \|\nabla d(\lag_1) - \nabla d(\lag_2)\|\\
= &\; \| \A \nabla f^*(-\A^\tp \lag_1) - \A \nabla f^*(-\A^\tp \lag_2)\|\\
\leq &\; \|\A\| \norm{\nabla f^*(-\A^\tp \lag_1) - \nabla f^*(-\A^\tp \lag_2)}\\
\leq &\; \frac{\|\A\|}{\monotone} \|\A^\tp(\lag_1 - \lag_2)\|\\
\leq &\; \frac{\|\A\|^2}{\monotone} \|\lag_1 - \lag_2\|
\end{aligned}
\end{equation}
which implies that $d(\lag)$ is $(\norm{A}^2/\monotone)$-smooth for all $ \lag \in \dom d $. Note that the dual function $ d $ is always differentiable given that $f$ is strongly convex.

To show that a differentiable function $f$ is $\monotone$-strongly convex, it is equivalent to show its gradient is strongly monotone, i.e., the inequality
\begin{equation*}
\langle \nabla f(\x) - \nabla f(\y), \x - \y \rangle \geq
\monotone \|\x - \y\|_2^2
\end{equation*}
holds for all $ \x,\y \in \dom f $.
We use this property to show strong convexity of $d(\lag)$.
The domain of $ d $ is characterized by \cref{eq:domain}.
For any $ \lag_1, \lag_2 \in \dom d $ and $\lag_1 \neq \lag_2$, we have
\begin{equation}\label{eq:monotoniciy}
\begin{aligned}
& \left\langle \nabla d(\lag_1) - \nabla d(\lag_2),%
\lag_1 - \lag_2 \right\rangle\\
= & \left\langle -\A\nabla f^*(-\A^\tp\lag_1)
+ \A\nabla f^*(-\A^\tp\lag_2), \lag_1 - \lag_2 \right\rangle\\
= & \left\langle \nabla f^*(-\A^\tp\lag_1) - \nabla f^*(-\A^\tp\lag_2),
- \A^\tp(\lag_1 - \lag_2) \right\rangle\\
\geq & \frac{1}{\smooth} \norm{\A^\tp(\lag_1 - \lag_2)}_2^2
\end{aligned}
\end{equation}
where the inequality is due to the strong convexity of $ f^* $.
According to~\eqref{eq:key-relation}, $\lag_1, \; \lag_2 \in\MC{R}(\A - \B \E^{-1} \C^\tp)$, then $\lag_1 - \lag_2 \in\MC{R}(\A - \B \E^{-1} \C^\tp)$.
There exists some $\VEC{u} \neq \VEC{0}$ such that
\begin{equation*}
\lag_1 - \lag_2 = (\A - \B \E^{-1} \C^\tp)\VEC{u}.
\end{equation*}
Using this relation and considering \( \A^\tp(\A - \B\E^{-1}\C^\tp) = \DCEC = \Lg \), \cref{eq:monotoniciy} reduces to
\begin{equation}\label{eq:monotone-bound}
\left\langle \nabla d(\lag_1) - \nabla d(\lag_2),%
\lag_1 - \lag_2 \right\rangle \geq%
\frac{1}{\smooth} \norm{\Lg \VEC{u}}^2.
\end{equation}
For a connected graph, the Laplacian always has an 1-dimensional null space spanned by $\VEC{1}$.
Thus, we can decompose $\VEC{u} = \VEC{u}_1 + \VEC{u}^\perp$ where $\VEC{u}_1 \in\mathbf{span}\{\VEC{1}\}$ and $\VEC{u}^\perp \perp \mathbf{span}\{\VEC{1}\}$.
First consider the case $\VEC{u}^\perp \neq \VEC{0}$.
We can further bound \cref{eq:monotoniciy} by
\begin{equation}\label{eq:bound-low}
\norm{\Lg \VEC{u}}^2
\geq \eigminnz^2(\Lg) \norm{\VEC{u}^\perp}^2.
\end{equation}
According to \cref{eq:key-relation}, we also have
\begin{multline*}
\|\lag_1 - \lag_2\|^2 = \left\| (\A - \B\E^{-1}\C^\tp) \VEC{u} \right\|^2
\leq \eigmax(\Lg) \left\| \VEC{u}^\perp \right\|^2
\end{multline*}
which immediately yields
\begin{equation}\label{eq:bound-max}
\norm{\VEC{u}^\perp}^2 \geq \frac{\|\lag_1 - \lag_2\|^2}{\eigmax(\Lg)}.
\end{equation}
We have used the fact that 
\begin{equation}
(\A - \B\E^{-1}\C^\tp)^\tp (\A - \B\E^{-1}\C^\tp) = \Lg.
\end{equation}
Combining \cref{eq:bound-low} and \cref{eq:bound-max}, we arrive at
\begin{equation}\label{eq:bound-inequality}
\left\lVert\Lg\VEC{u} \right\rVert^2
\geq \frac{\eigmin^2(\Lg)}{\eigmax(\Lg)} \|\lag_1 - \lag_2\|^2
\end{equation}
which, together with \cref{eq:monotone-bound}, yields
\begin{equation}\label{eq:monotone-result}
\left\langle \nabla d(\lag_1) - \nabla d(\lag_2),%
\lag_1 - \lag_2 \right\rangle \geq%
\frac{1}{\smooth} \frac{\eigmin^2(\Lg)}{\eigmax(\Lg)} \|\lag_1 - \lag_2\|^2.
\end{equation}
Now consider the case \( \VEC{u}^\perp = \VEC{0} \).
Since
\begin{equation*}
\norm{\lag_1 - \lag_2}^2 = \norm{(\A - \B\E^{-1}\C^\tp) \VEC{u}}^2 = 0,
\end{equation*}
we have $\lag_1 = \lag_2$ and \cref{eq:monotone-result} holds.
As a result, we have
\begin{multline}
\langle \nabla d(\lag_1) - \nabla d(\lag_2), \lag_1 - 
\lag_2 \rangle
\geq \frac{1}{\smooth} \frac{\eigmin^2(\Lg)}{\eigmax(\Lg)} \|\lag_1 - \lag_2\|^2
\end{multline}
for all $\lag_1, \; \lag_2 \in \dom d$.
Therefore, $d$ is $\hat{\alpha}$-strongly convex, with $\hat{\alpha}$ defined by
\begin{equation*}
\hat{\alpha} = \frac{1}{\smooth} \frac{\eigminnz^2(\Lg)}{\eigmax(\Lg)}
\end{equation*}
The proof is completed.

\section{Proof of Theorem~\ref{th:main-result}}%
\label{sec:proof-main-result}
	Similar to the proof of DR algorithm in the primal domain, we first show that $ \refl{d_1} $ is a contraction.
	This is due to \cref{th:dual-properties} and \cref{th:contraction}, and the contraction factor is
	\begin{equation}
	\dualform{c} = \max \left\{ \frac{\gamma \dualform{\smooth} - 1}{\gamma \dualform{\smooth} + 1}, \frac{1 - \gamma\dualform{\monotone}}{1 + \gamma\dualform{\monotone}} \right\}
	\end{equation}
	where $\dualform{\smooth}, \dualform{\monotone}$ are defined in \cref{th:dual-properties}.
	For any $ x, y \in \dom\: d_1 $, $ \|\refl{d_1} x - \refl{d_1} y\| \leq \dualform{c} \| x - y\| $.
	Given that $d_2 = g_1^* = g_2$ is closed, convex and proper, then $\refl{d_2}$ is nonexpansive.
	Hence, the composition $\refl{d_2}\refl{d_1}$ is also contractive with contraction factor $\dualform{c}$.
	Let $\opT = (1-\average)\opI + \average\refl{d_1}\refl{d_2}$ be the Dougals-Rachford operator and $\bar{\dr}$ a fixed point of $\opT$. Thus,
	\begin{equation*}
	\begin{aligned}
	& \|\dr^{k+1} - \dualform{\dr}\| = \|\opT \dr^k - \opT \bar{\dr}\|\\
	= & \|(1-\average) ( \dr^k - \bar{\dr}) + \average ( \refl{d_1} \refl{d_2} \dr^k - \refl{d_1} \refl{d_2} \bar{\dr} ) \|\\
	\leq & |1-\average| \|\dr^k - \dualform{\dr}\| + \average \| \refl{d_1} \refl{d_2} \dr^k - \refl{d_1} \refl{d_2} \bar{\dr} \|\\
	\leq & (|1 - \average| + \average \dualform{c}) \|\dr^k - \dualform{\dr}\|.
	\end{aligned}
	\end{equation*}
	Furthermore, it is easy to see that 
	\begin{equation*}
	|1 - \average| + \average \dualform{c} < 1
	\end{equation*}
	if and only if $ \average \in (0, 2/(1+\dualform{c})) $.

\section{Proof of Corollary~\ref{th:dual-property-equivalence}}%
\label{sec:proof-equivalence}

The structure of matrix $\A$ implies that full row rank is attainable only when all nodes can be assigned into a single group, i.e., there exists a \emph{central node} connecting all others.
In such cases, we have $\A = \Id$, $\B = \VEC{1}$.
As a consequence, \cref{th:matrix-relation} leads to $\C = \Id$, $\D = \Id$ and $\E = n$.
The group degree matrix $\E$ reduces to a scalar since there is only one group.
The Laplacian of the communication graph can be obtained consequently
\begin{equation*}
\Lg = \DCEC = \Id - \frac{1}{n} \VEC{1}\VEC{1}^\tp
\end{equation*}
which happens to be the Laplacian matrix of a complete graph with $n$ nodes, $K_n$, up to a scale, i.e., $ \mat{L}_{K_n} = n \Lg $.
Since $\mat{L}_{K_n}$ has all only one nonzero eigenvalues $n$ with multiplicity $n-1$, i.e., $\eigminnz(\mat{L}_{K_n}) = \eigmax(\mat{L}_{K_{n}})$. 
Therefore, we have $\eigminnz(\Lg) = \eigmax(\Lg) = 1$, and we arrive at
\begin{equation*}
\dualform{\monotone} = \frac{\eigminnz^2(\Lg)}{\smooth\eigmax(\Lg)}
= \frac{1}{\smooth}.
\end{equation*}
Considering $\sigmin(\A) = 1$, we have completed the proof.


\begin{thebibliography}{10}
\providecommand{\url}[1]{#1}
\csname url@samestyle\endcsname
\providecommand{\newblock}{\relax}
\providecommand{\bibinfo}[2]{#2}
\providecommand{\BIBentrySTDinterwordspacing}{\spaceskip=0pt\relax}
\providecommand{\BIBentryALTinterwordstretchfactor}{4}
\providecommand{\BIBentryALTinterwordspacing}{\spaceskip=\fontdimen2\font plus
\BIBentryALTinterwordstretchfactor\fontdimen3\font minus
  \fontdimen4\font\relax}
\providecommand{\BIBforeignlanguage}[2]{{%
\expandafter\ifx\csname l@#1\endcsname\relax
\typeout{** WARNING: IEEEtran.bst: No hyphenation pattern has been}%
\typeout{** loaded for the language `#1'. Using the pattern for}%
\typeout{** the default language instead.}%
\else
\language=\csname l@#1\endcsname
\fi
#2}}
\providecommand{\BIBdecl}{\relax}
\BIBdecl

\bibitem{nedic2009distributed}
A.~Nedic and A.~Ozdaglar, ``Distributed {{Subgradient Methods}} for
  {{Multi}}-{{Agent Optimization}},'' \emph{IEEE Transactions on Automatic
  Control}, vol.~54, no.~1, pp. 48--61, Jan. 2009.

\bibitem{nedic2010constrained}
A.~Nedic, A.~Ozdaglar, and P.~A. Parrilo, ``Constrained consensus and
  optimization in multi-agent networks,'' \emph{IEEE Transactions on Automatic
  Control}, vol.~55, no.~4, pp. 922--938, 2010.

\bibitem{schizas2008consensus}
I.~D. Schizas, A.~Ribeiro, and G.~B. Giannakis, ``Consensus in {{Ad Hoc WSNs
  With Noisy Links}} -- {{Part I}}: {{Distributed Estimation}} of
  {{Deterministic Signals}},'' \emph{IEEE Transactions on Signal Processing},
  vol.~56, no.~1, pp. 350--364, Jan. 2008.

\bibitem{ren2007information}
W.~Ren, R.~W. Beard, and E.~M. Atkins,
  ``\BIBforeignlanguage{English}{Information consensus in multivehicle
  cooperative control},'' \emph{\BIBforeignlanguage{English}{IEEE Control Syst.
  Mag.}}, vol.~27, no.~2, pp. 71--82, Apr. 2007.

\bibitem{boyd2011distributed}
S.~Boyd, N.~Parikh, E.~Chu, B.~Peleato, and J.~Eckstein, ``Distributed
  {{Optimization}} and {{Statistical Learning}} via the {{Alternating Mirection
  Method}} of {{Multipliers}},'' \emph{Foundations and Trends\textregistered{}
  in Machine Learning}, vol.~3, no.~1, pp. 1--122, 2011.

\bibitem{giannakis2016decentralized}
G.~B. Giannakis, Q.~Ling, G.~Mateos, I.~D. Schizas, and H.~Zhu, ``Decentralized
  {{Learning}} for {{Wireless Communications}} and {{Networking}},'' in
  \emph{Splitting {{Methods}} in {{Communication}}, {{Imaging}}, {{Science}},
  and {{Engineering}}}, R.~Glowinski, S.~J. Osher, and W.~Yin, Eds.\hskip 1em
  plus 0.5em minus 0.4em\relax Cham: {Springer}, 2016, pp. 461--497.

\bibitem{olfati-saber2007consensus}
R.~{Olfati-Saber}, J.~A. Fax, and R.~M. Murray, ``Consensus and {{Cooperation}}
  in {{Networked Multi}}-{{Agent Systems}},'' \emph{Proceedings of the IEEE},
  vol.~95, no.~1, pp. 215--233, Jan. 2007.

\bibitem{sayed2014adaptation}
A.~H. Sayed, ``\BIBforeignlanguage{English}{Adaptation, {{Learning}}, and
  {{Optimization}} over {{Networks}}},''
  \emph{\BIBforeignlanguage{English}{MAL}}, vol.~7, no. 4-5, pp. 311--801, Jul.
  2014.

\bibitem{rabbat2004distributed}
M.~Rabbat and R.~Nowak, ``\BIBforeignlanguage{en}{Distributed optimization in
  sensor networks},'' in \emph{\BIBforeignlanguage{en}{Proceedings of the Third
  International Symposium on {{Information}} Processing in Sensor Networks -
  {{IPSN}}'04}}.\hskip 1em plus 0.5em minus 0.4em\relax Berkeley, California,
  USA: {ACM Press}, 2004, p.~20.

\bibitem{sundharram2010distributed}
S.~Sundhar~Ram, A.~Nedi\'c, and V.~V. Veeravalli,
  ``\BIBforeignlanguage{en}{Distributed {{Stochastic Subgradient Projection
  Algorithms}} for {{Convex Optimization}}},'' \emph{\BIBforeignlanguage{en}{J
  Optim Theory Appl}}, vol. 147, no.~3, pp. 516--545, Dec. 2010.

\bibitem{duchi2012dual}
J.~C. Duchi, A.~Agarwal, and M.~J. Wainwright, ``Dual {{Averaging}} for
  {{Distributed Optimization}}: {{Convergence Analysis}} and {{Network
  Scaling}},'' \emph{IEEE Transactions on Automatic Control}, vol.~57, no.~3,
  pp. 592--606, Mar. 2012.

\bibitem{dimakis2010gossip}
A.~G. Dimakis, S.~Kar, J.~M.~F. Moura, M.~G. Rabbat, and A.~Scaglione, ``Gossip
  {{Algorithms}} for {{Distributed Signal Processing}},'' \emph{Proceedings of
  the IEEE}, vol.~98, no.~11, pp. 1847--1864, Nov. 2010.

\bibitem{bertsekas1989parallel}
D.~P. Bertsekas and J.~N. Tsitsiklis, \emph{Parallel and Distributed
  Computation: Numerical Methods}.\hskip 1em plus 0.5em minus 0.4em\relax
  Englewood Cliffs, NJ: {Prentice-Hall}, 1989, vol.~23.

\bibitem{mateos2009distributed}
G.~Mateos, I.~D. Schizas, and G.~B. Giannakis, ``Distributed {{Recursive
  Least}}-{{Squares}} for {{Consensus}}-{{Based In}}-{{Network Adaptive
  Estimation}},'' \emph{IEEE Transactions on Signal Processing}, vol.~57,
  no.~11, pp. 4583--4588, Nov. 2009.

\bibitem{ma2018hybrid}
M.~Ma, A.~N. Nikolakopoulos, and G.~B. Giannakis,
  ``\BIBforeignlanguage{en}{Hybrid {{ADMM}}: A unifying and fast approach to
  decentralized optimization},'' \emph{\BIBforeignlanguage{en}{EURASIP J. Adv.
  Signal Process.}}, vol. 2018, no.~1, p.~73, Dec. 2018.

\bibitem{glowinski1975approximation}
R.~Glowinski and A.~Marroco, ``Sur l'approximation, par \'el\'ements finis
  d'ordre un, et la r\'esolution, par p\'enalisation-dualit\'e d'une classe de
  probl\`emes de {{Dirichlet}} non lin\'eaires,'' \emph{Revue fran{\c c}aise
  d'automatique, informatique, recherche op\'erationnelle. Analyse
  num\'erique}, vol.~9, no.~R2, pp. 41--76, 1975.

\bibitem{iutzeler2016explicit}
F.~Iutzeler, P.~Bianchi, P.~Ciblat, and W.~Hachem, ``Explicit {{Convergence
  Rate}} of a {{Distributed Alternating Direction Method}} of
  {{Multipliers}},'' \emph{IEEE Transactions on Automatic Control}, vol.~61,
  no.~4, pp. 892--904, Apr. 2016.

\bibitem{gabay1976dual}
D.~Gabay and B.~Mercier, ``A dual algorithm for the solution of nonlinear
  variational problems via finite element approximation,'' \emph{Computers \&
  Mathematics with Applications}, vol.~2, no.~1, pp. 17--40, Jan. 1976.

\bibitem{hong2017linear}
M.~Hong and Z.-Q. Luo, ``\BIBforeignlanguage{en}{On the linear convergence of
  the alternating direction method of multipliers},''
  \emph{\BIBforeignlanguage{en}{Math. Program.}}, vol. 162, no. 1-2, pp.
  165--199, Mar. 2017.

\bibitem{deng2015global}
W.~Deng and W.~Yin, ``\BIBforeignlanguage{en}{On the {{Global}} and {{Linear
  Convergence}} of the {{Generalized Alternating Direction Method}} of
  {{Multipliers}}},'' \emph{\BIBforeignlanguage{en}{J Sci Comput}}, vol.~66,
  no.~3, pp. 889--916, May 2015.

\bibitem{shi2014linear}
W.~Shi, Q.~Ling, K.~Yuan, G.~Wu, and W.~Yin, ``On the {{Linear Convergence}} of
  the {{ADMM}} in {{Decentralized Consensus Optimization}},'' \emph{IEEE
  Transactions on Signal Processing}, vol.~62, no.~7, pp. 1750--1761, Apr.
  2014.

\bibitem{makhdoumi2017convergence}
A.~Makhdoumi and A.~Ozdaglar, ``Convergence {{Rate}} of {{Distributed ADMM Over
  Networks}},'' \emph{IEEE Transactions on Automatic Control}, vol.~62, no.~10,
  pp. 5082--5095, Oct. 2017.

\bibitem{lin2015global}
T.~Lin, S.~Ma, and S.~Zhang, ``On the {{Global Linear Convergence}} of the
  {{ADMM}} with {{MultiBlock Variables}},'' \emph{SIAM J. Optim.}, vol.~25,
  no.~3, pp. 1478--1497, Jan. 2015.

\bibitem{franca2017how}
G.~Fran{\c c}a and J.~Bento, ``How is {{Distributed ADMM Affected}} by
  {{Network Topology}}?'' \emph{arXiv:1710.00889 [math, stat]}, Oct. 2017.

\bibitem{ghadimi2014admm}
E.~Ghadimi, A.~Teixeira, M.~G. Rabbat, and M.~Johansson, ``The {{ADMM}}
  algorithm for distributed averaging: {{Convergence}} rates and optimal
  parameter selection,'' in \emph{Proc. of {{Asilomar Conf}}. on {{Signals}},
  {{Systems}}, and {{Computers}}}, Nov. 2014, pp. 783--787.

\bibitem{ling2016weighted}
Q.~Ling, Y.~Liu, W.~Shi, and Z.~Tian, ``Weighted {{ADMM}} for {{Fast
  Decentralized Network Optimization}},'' \emph{IEEE Transactions on Signal
  Processing}, vol.~64, no.~22, pp. 5930--5942, Nov. 2016.

\bibitem{ma2018graphaware}
M.~Ma and G.~B. Giannakis, ``Graph-aware {{Weighted Hybrid ADMM}} for {{Fast
  Decentralized Optimization}},'' in \emph{2018 52nd {{Asilomar Conference}} on
  {{Signals}}, {{Systems}}, and {{Computers}}}, Oct. 2018, pp. 1881--1885.

\bibitem{rockafellar1976monotone}
R.~Rockafellar, ``Monotone {{Operators}} and the {{Proximal Point
  Algorithm}},'' \emph{SIAM J. Control Optim.}, vol.~14, no.~5, pp. 877--898,
  Aug. 1976.

\bibitem{bauschke2017convex}
H.~H. Bauschke and P.~L. Combettes, \emph{Convex {{Analysis}} and {{Monotone
  Operator Theory}} in {{Hilbert Spaces}}}, 2nd~ed., ser. {{CMS Books}} in
  {{Mathematics}}.\hskip 1em plus 0.5em minus 0.4em\relax Cham: {Springer
  International Publishing}, 2017.

\bibitem{lions1979splitting}
P.~Lions and B.~Mercier, ``Splitting {{Algorithms}} for the {{Sum}} of {{Two
  Nonlinear Operators}},'' \emph{SIAM J. Numer. Anal.}, vol.~16, no.~6, pp.
  964--979, Dec. 1979.

\bibitem{eckstein1989splitting}
J.~Eckstein, ``\BIBforeignlanguage{eng}{Splitting methods for monotone
  operators with applications to parallel optimization},'' Thesis,
  Massachusetts Institute of Technology, 1989.

\bibitem{eckstein1992douglasrachford}
J.~Eckstein and D.~P. Bertsekas, ``\BIBforeignlanguage{en}{On the
  {{Douglas}}-{{Rachford}} splitting method and the proximal point algorithm
  for maximal monotone operators},'' \emph{\BIBforeignlanguage{en}{Mathematical
  Programming}}, vol.~55, no. 1-3, pp. 293--318, Apr. 1992.

\bibitem{douglas1956numerical}
J.~Douglas and H.~H. Rachford, ``On the {{Numerical Solution}} of {{Heat
  Conduction Problems}} in {{Two}} and {{Three Space Variables}},''
  \emph{Transactions of the American Mathematical Society}, vol.~82, no.~2, pp.
  421--439, 1956.

\bibitem{gabay1983applications}
D.~Gabay, ``Applications of the {{Method}} of {{Multipliers}} to {{Variational
  Inequalities}},'' in \emph{Studies in {{Mathematics}} and {{Its
  Applications}}}, ser. Augmented {{Lagrangian Methods}}: {{Applications}} to
  the {{Numerical Solution}} of {{Boundary}}-{{Value Problems}}, M.~Fortin and
  R.~Glowinski, Eds.\hskip 1em plus 0.5em minus 0.4em\relax Amsterdam:
  {North-Holland}, Jan. 1983, vol.~15, pp. 299--331.

\bibitem{giselsson2017linear}
P.~Giselsson and S.~Boyd, ``Linear {{Convergence}} and {{Metric Selection}} for
  {{Douglas}}-{{Rachford Splitting}} and {{ADMM}},'' \emph{IEEE Transactions on
  Automatic Control}, vol.~62, no.~2, pp. 532--544, Feburary 2017.

\bibitem{bolla1993spectra}
M.~Bolla, ``Spectra, {{Euclidean}} representations and clusterings of
  hypergraphs,'' \emph{Discrete Mathematics}, vol. 117, no.~1, pp. 19--39, Jul.
  1993.

\bibitem{ryu2016primer}
E.~K. Ryu and S.~Boyd, ``\BIBforeignlanguage{en}{A {{Primer}} on {{Monotone
  Operator Methods}}},'' \emph{\BIBforeignlanguage{en}{Appl. Comput. Math}},
  vol.~15, no.~1, pp. 3--43, 2016.

\bibitem{brandes2001faster}
U.~Brandes, ``A faster algorithm for betweenness centrality,'' \emph{The
  Journal of Mathematical Sociology}, vol.~25, no.~2, pp. 163--177, Jun. 2001.

\end{thebibliography}
\end{document}